\documentclass[12pt]{article}%
\usepackage{amsfonts}
\usepackage{amsmath,amssymb}
\usepackage[applemac]{inputenc}
\usepackage{amsmath,amssymb,fullpage}
\usepackage{showlabels}
\usepackage{color}
\usepackage{amsmath}
\usepackage{amssymb}
\usepackage{graphicx}%
\providecommand{\U}[1]{\protect\rule{.1in}{.1in}}

\newtheorem{definition}{Definition}[section]
\newtheorem{example}{Example}[section]
\newtheorem{theorem}[definition]{Theorem}

\newtheorem{remark}[definition]{ \it Remark}

\newtheorem{proposition}[definition]{Proposition}
\newtheorem{lemma}[definition]{Lemma}

\numberwithin{equation}{section}

\def\1B{\text{1\!\!I}}

\begin{document}

\title{Stochastic optimal control of McKean-Vlasov equations with anticipating law}
\author{\ Nacira A\textsc{GRAM}\thanks{Department of Mathematics, University of Oslo,
P.O. Box 1053 Blindern, N--0316 Oslo, Norway. Email:
\texttt{naciraa@math.uio.no.}} \thanks{University Mohamed Khider of Biskra,
Algeria.} \thanks{This research was carried out with support of the Norwegian
Research Council, within the research project Challenges in Stochastic
Control, Information and Applications (STOCONINF), project number 250768/F20.}}
\date{30 June 2017}
\maketitle

\begin{abstract}
We are interested in Pontryagin's stochastic maximum principle of controlled
McKean-Vlasov stochastic differential equations (SDE). We allow the law to be
anticipating, in the sense that, the coefficients (the drift and the diffusion
coefficient) depend not only of the solution at the current time $t$, but also
on the law of the future values of the solution $P_{X(t+\delta)}$, for a given
positive constant $\delta$. We emphasise that being anticipating w.r.t. the
law of the solution process does not mean being anticipative in the sense that
it anticipates the driving Brownian motion. As an adjoint equation, a new type
of delayed backward stochastic differential equations (BSDE) with implicit
terminal condition is obtained.

By using that the expectation of any random variable is a function of its law,
our BSDE can be written in a simple form. Then, we prove existence and
uniqueness of the solution of the delayed BSDE with implicit terminal value,
i.e. with terminal value being a function of the law of the solution itself.

\end{abstract}


\section{Introduction}

The stochastic maximum principle and the characterization of the optimal
control by Pontryagin's maximum principle have been studied for the classical
case by many authors such as Bismut \cite{Bismut}, and Bensoussan
\cite{Bensoussan1}, \cite{Bensoussan2} and so on. Recently, this principle was
extended to more general cases such as controlled MacKean-Vlasov systems by
Carmona and Delarue \cite{carmona1}, \cite{carmona}, in the sense that both
the drift and the diffusion coefficients in their dynamics, are supposed to
depend at time $t$ on the solution $X(t)$ and on its law $P_{X(t)}$ as follows%
\[
\left\{
\begin{array}
[c]{ll}%
dX(t) & =b(t,X(t),P_{X(t)},u(t))dt+\sigma(t,X(t),P_{X(t)},u(t))dB(t),t\in
\lbrack0,T],\\
X(0) & =x.
\end{array}
\right.
\]
In this paper, we want to extend Pontryagin's stochastic maximum principle to
a more general case, the case of controlled McKean-Vlasov SDE with
anticipating law, i.e., our dynamics are assumed to satisfy for a given
positive constant $\delta$, the equation%
\begin{equation}
\left\{
\begin{array}
[c]{ll}%
dX(t) & =b(t,X(t),P_{X(t+\delta)},u(t))dt+\sigma(t,X(t),P_{X(t+\delta
)},u(t))dB(t),t\in\lbrack0,T],\\
X(0) & =x,\\
X(t) & =X(T);t\geq T,
\end{array}
\right.  \label{*}%
\end{equation}
where $u(t)$ is our control process. Here the coefficients depend at time $t$
on both the solution $X(t)$ and the law $P_{X(t+\delta)}$ of the anticipating
solution. We remark here that the SDE $\left(  \ref{*}\right)  $ being
anticipative w.r.t. the law of the solution process does not mean being
anticipative in the sense that it anticipates the driving Brownian motion $B$.

The performance functional for a control $u$ is given by
\[
J\left(  u\right)  =\mathbb{E[}g\left(  X\left(  T\right)  ,P_{X(T)}\right)  +%
{\textstyle\int_{0}^{T}}
l\left(  s,X\left(  s\right)  ,P_{X(s+\delta)},u\left(  s\right)  \right)
ds]\text{,}%
\]
for some given bounded functions $g$ and $l$ and we want to maximize this
performance over the set $\mathcal{U}$ of admissible control processes (will
be specified later), as follows: Find $u^{\ast}\in\mathcal{U}$ such that
$J(u^{\ast})=\sup_{u\in\mathcal{U}}J\left(  u\right)  $.

We define the Hamiltonian $H$ associated to this problem to be%
\[
H(x,\mu,u,p,q)=l(t,x,\mu,u)+b(t,x,\mu,u)p+\sigma(t,x,\mu,u)q,
\]
and we show that the couple $(p,q)$ is the solution of the adjoint backward
stochastic differential equation given by%
\begin{equation}%
\begin{array}
[c]{l}%
dp(t)=-\{(\partial_{x}b)(t,X^{\ast}\left(  t\right)  ,P_{X^{^{\ast}}\left(
t+\delta\right)  },u^{\ast}\left(  t\right)  )p(t)\\
+(\partial_{x}\sigma)(t,X^{\ast}\left(  t\right)  ,P_{X^{^{\ast}}\left(
t+\delta\right)  },u^{\ast}\left(  t\right)  )q(t)+\left(  \partial
_{x}l\right)  (t,X^{\ast}(t),P_{X^{\ast}(t+\delta)},u^{\ast}(t))\\
+\mathbb{\tilde{E}[(}\partial_{\mu}b)(t-\delta,\tilde{X}^{\ast}\left(
t-\delta\right)  ,P_{X^{\ast}\left(  t\right)  },X^{\ast}\left(  t\right)
,\tilde{u}^{\ast}\left(  t-\delta\right)  )\tilde{p}\left(  t-\delta\right)
]I_{\left[  \delta,T\right]  }\left(  t\right) \\
+\mathbb{\tilde{E}[(}\partial_{\mu}\sigma)(t-\delta,\tilde{X}^{\ast}\left(
t-\delta\right)  ,P_{X^{\ast}\left(  t\right)  },X^{\ast}\left(  t\right)
,\tilde{u}^{\ast}\left(  t-\delta\right)  )\tilde{q}\left(  t-\delta\right)
]I_{\left[  \delta,T\right]  }\left(  t\right) \\
+\mathbb{\tilde{E}[(}\partial_{\mu}l)(t-\delta,\tilde{X}^{\ast}\left(
t-\delta\right)  ,P_{X^{\ast}\left(  t\right)  },X^{\ast}\left(  t\right)
,\tilde{u}^{\ast}\left(  t-\delta\right)  )]I_{\left[  \delta,T\right]
}\left(  t\right)  \}dt+q(t)dB(t)\text{, }t\in\left[  0,T\right]  \text{,}%
\end{array}
\label{n1}%
\end{equation}
with terminal condition%
\begin{equation}%
\begin{array}
[c]{ll}%
p(T)= & (\partial_{x}g)(X^{\ast}(T),P_{X^{\ast}(T)})+\mathbb{\tilde{E}%
[(}\partial_{\mu}g)(\tilde{X}^{\ast}(T),P_{X^{\ast}(T)},X^{\ast}(T))]\\
& +%
{\textstyle\int_{T-\delta}^{T}}
(\mathbb{\tilde{E}[(}\partial_{\mu}b)(t,\tilde{X}^{\ast}\left(  t\right)
,P_{X^{\ast}\left(  T\right)  },X^{\ast}\left(  T\right)  ,\tilde{u}^{\ast
}\left(  t\right)  )\tilde{p}\left(  t\right)  ]\\
& +\mathbb{\tilde{E}[(}\partial_{\mu}\sigma)(t,\tilde{X}^{\ast}\left(
t\right)  ,P_{X^{\ast}\left(  T\right)  },X^{\ast}\left(  T\right)  ,\tilde
{u}^{\ast}\left(  t\right)  )\tilde{q}\left(  t\right)  ]\\
& +\mathbb{\tilde{E}[(}\partial_{\mu}l)(t,\tilde{X}^{\ast}(t),P_{X^{\ast}%
(T)},X^{\ast}(T),\tilde{u}^{\ast}(t))])dt.
\end{array}
\label{n2}%
\end{equation}
The adjoint equation $\left(  \ref{n1}\right)  -\left(  \ref{n2}\right)  $ is
a new type of delayed McKean-Vlasov BSDE with implicit terminal
condition,i.e., with terminal value being a function of the law of the
solution itself. In order to write it in a more comprehensible form, we use
the fact that the expectation of any random variable is a function of its law,
and under suitable assumptions on both the driver and a terminal value, we can
get existence and uniqueness of our delayed BSDE. It is a generalisation of
the adjoint equation of the above mentioned problem we are interested
in.\newline

Stochastic Pontryagin's maximum principle in both cases partial and complete
information of mean-field systems has been studied for example by Anderson and
Djehiche in \cite{AD} and Hu el al in \cite{HOS} and for more details about
mean-field systems, we refer to Lions \cite{Lions}, Cardaliaguet notes
\cite{Cardaliaguet} and Buckdahn et al \cite{BLPR}.\newline Delayed BSDEs have
been studied by Delong and Imkeller \cite{delond}; they have later been
extended by the same authors to the jump case, and studied \cite{del} by the
help of the Malliavin calculus and for more details, we refer to Delong's book
\cite{de}.\newline

To the best of our knowledge our paper is the first to study optimal control
problems of mean-field SDEs with anticipating law.

The paper is organized as follows: In the next section, we give some
preliminaries which will be used throughout this work. In section $3$, the
existence and the uniqueness of McKean-Vlasov SDEs with anticipating law is
investigated. Section $4$ is devoted to the study of Pontryagin's stochastic
maximum principle. In the last section, we prove the existence and the
uniqueness for the associated delayed McKean-Vlasov BSDEs with implicit
terminal condition.\newline\newline

\textit{This work has been presented at seminars and conferences in Brest,
Biskra, Marrakech, Mans and Oslo.}

\section{Framework}

We introduce some notations, definitions and spaces which will be used
throughout this work. Let $\left(  \Omega,\mathcal{F},P\right)  $ be a
complete probability space, $B$ a $d$-dimensional Brownian motion and
$\mathbb{F}=\left(  \mathcal{F}_{t}\right)  _{t\geq0}$ the Brownian filtration
generated by $B$ and completed by all $P$-null sets$.$ Let $\mathcal{P}%
_{2}(\mathbb{R}^{d}):=\{\mu\in\mathcal{P}(\mathbb{R}^{d}):%
{\textstyle\int_{\mathbb{R}^{d}}}
|x|^{2}\mu(dx)<+\infty\},$ where $\mathcal{P}(\mathbb{R}^{d})$ is the space of
all the probability measures on $(\mathbb{R}^{d},\mathcal{B}(\mathbb{R}^{d}%
))$; recall that $\mathcal{B}(\mathbb{R}^{d})$ denotes the Borel $\sigma
$-field over $\mathbb{R}^{d}$. We endow $\mathcal{P}_{2}(\mathbb{R}^{d})$ with
the $2$-Wasserstein metric $W_{2}$ on $\mathcal{P}_{2}(\mathbb{R}^{d})$: For
$\mu_{1},\mu_{2}\in\mathcal{P}_{2}(\mathbb{R}^{d})$, the $2$-Wasserstein
distance is defined by
\[%
\begin{array}
[c]{l}%
W_{2}\left(  \mu_{1},\mu_{2}\right)  =\inf\{(%
{\textstyle\int_{\mathbb{R}^{d}}}
\left\vert x-y\right\vert ^{2}\mu(dx,dy))^{\frac{1}{2}}:\mu\in\mathcal{P}%
_{2}(\mathbb{R}^{d}\times\mathbb{R}^{d})\\
\text{ \ \ \ \ \ \ \ \ \ \ \ \ \ \ \ \ \ \ \ \ \ \ \ \ \ \ \ \ \ \ with }%
\mu(\mathbb{\cdot}\times\mathbb{R}^{d}):=\mu_{1},\text{ }\mu(\mathbb{R}%
^{d}\times\mathbb{\cdot}):=\mu_{2}\}\text{.}%
\end{array}
\]
We also remark that, if $\left(  \Omega,\mathcal{F},P\right)  $ is "rich
enough" in the sense that
\[
\mathcal{P}_{2}(\mathbb{R}^{d}\times\mathbb{R}^{d})=\left\{  P_{\zeta}\text{,
}\zeta\in L^{2}\left(  \Omega,\mathcal{F},P;\mathbb{R}^{d}\times\mathbb{R}%
^{d}\right)  \right\}  \text{,}%
\]
then we also have
\[
W_{2}\left(  \mu_{1},\mu_{2}\right)  =\inf\{(\mathbb{E[}\left\vert \zeta
-\eta\right\vert ^{2}])^{\frac{1}{2}}\text{, }\zeta\text{, }\eta\in
L^{2}\left(  \Omega,\mathcal{F},P;\mathbb{R}^{d}\right)  \text{,with }%
P_{\zeta}=\mu_{1}\text{, }P_{\eta}=\mu_{2}\}\text{.}%
\]
Let $\left(  \tilde{\Omega},\mathcal{\tilde{F}},\tilde{P}\right)  :=\left(
\Omega,\mathcal{F},P\right)  $, $\mathbb{\tilde{F}}:=\mathbb{F}$ and $\left(
\bar{\Omega},\mathcal{\bar{F}},\bar{P}\right)  =\left(  \Omega,\mathcal{F}%
,P\right)  \otimes\left(  \tilde{\Omega},\mathcal{\tilde{F}},\tilde{P}\right)
$. For any measurable space $\left(  E,\mathcal{E}\right)  $ and any random
variable $\zeta:$ $\left(  \Omega,\mathcal{F},P\right)  \rightarrow\left(
E,\mathcal{E}\right)  $, we put $\tilde{\zeta}\left(  \tilde{\omega}\right)
:=\zeta\left(  \tilde{\omega}\right)  ,$ $\tilde{\omega}\in\tilde{\Omega
}=\Omega$, $\zeta\left(  \omega,\tilde{\omega}\right)  :=\zeta\left(
\omega\right)  $, $\tilde{\zeta}\left(  \omega,\tilde{\omega}\right)
:=\tilde{\zeta}\left(  \tilde{\omega}\right)  $, $\left(  \omega,\tilde
{\omega}\right)  \in\Omega\times\tilde{\Omega}.$ We observe that $\tilde
{\zeta}$ on $\left(  \tilde{\Omega},\mathcal{\tilde{F}},\tilde{P}\right)  $ is
a copy of $\zeta$ on $\left(  \Omega,\mathcal{F},P\right)  ,$ and$\ \zeta$,
$\tilde{\zeta}$ are i.i.d under $\bar{P}.$ Moreover, for $\zeta$,
$\eta:\left(  \Omega,\mathcal{F},P\right)  \rightarrow\left(  E,\mathcal{E}%
\right)  $ random variables and $\varphi:\left(  E^{2},\mathcal{E}^{2}\right)
\rightarrow\left(  B,\mathcal{B}(%
\mathbb{R}
)\right)  $ a bounded and measurable function, we have%
\begin{align*}
\mathbb{\tilde{E}[}\varphi(\tilde{\zeta},\eta)]\left(  \omega\right)   &  =%
{\textstyle\int_{\tilde{\Omega}}}
\varphi(\tilde{\zeta}\left(  \omega\right)  ,\eta\left(  \omega\right)
)\tilde{P}\left(  d\tilde{\omega}\right)  \\
&  =\mathbb{E}\left[  \varphi\left(  \zeta,y\right)  \right]  _{\diagup
y=\eta\left(  \omega\right)  }.
\end{align*}
We recall now the notion of derivative of a function $\varphi:\mathcal{P}%
_{2}(\mathbb{R}^{d})\rightarrow\mathbb{R}$ w.r.t a probability measure $\mu$,
which was studied by Lions in his course at Collège de France in \cite{Lions};
see also the notes of Cardaliaguet in \cite{Cardaliaguet}, the works by
Carmona and Delarue in \cite{carmona} and in Buckdahn, Li, Peng and Rainer
\cite{BLPR}. We say that $\varphi$ is differentiable at $\mu$ if, for the
lifted function $\tilde{\varphi}(\zeta):=\varphi(P_{\zeta})$, $\zeta\in
L^{2}\left(  \Omega,\mathcal{F},P;\mathbb{R}^{d}\right)  $, there is some
$\zeta_{0}\in L^{2}\left(  \Omega,\mathcal{F},P;\mathbb{R}^{d}\right)  $ with
$P_{\zeta_{0}}=\mu$, such that $\tilde{\varphi}$ is differentiable in the
Frèchet sense at $\zeta_{0}$, such that there exists a linear continuous
mapping $D\tilde{\varphi}(\zeta_{0}):L^{2}\left(  \Omega,\mathcal{F}%
,P;\mathbb{R}^{d}\right)  \rightarrow\mathbb{R}$ $(L\left(  L^{2}\left(
\Omega,\mathcal{F},P;\mathbb{R}^{d}\right)  ;\mathbb{R}\right)  )$, such that
\begin{align*}
\varphi(P_{\zeta_{0}+\eta})-\varphi(P_{\zeta_{0}}) &  =\tilde{\varphi}%
(\zeta_{0}+\eta)-\tilde{\varphi}(\zeta_{0})\\
&  =(D\tilde{\varphi})(\zeta_{0})(\eta)+o(\left\vert \eta\right\vert
_{L^{2}\left(  \Omega\right)  }^{2})\text{,}%
\end{align*}
for $\left\vert \eta\right\vert _{L^{2}\left(  \Omega\right)  }^{2}%
\rightarrow0$, $\eta\in L^{2}\left(  \Omega,\mathcal{F},P;\mathbb{R}%
^{d}\right)  .$ With the identification that \newline$L\left(  L^{2}\left(
\Omega,\mathcal{F},P;\mathbb{R}^{d}\right)  ;\mathbb{R}\right)  \equiv
L^{2}\left(  \Omega,\mathcal{F},P;\mathbb{R}^{d}\right)  $, given by Riesz'
representation theorem, we can write%
\[
\varphi(P_{\zeta_{0}+\eta})-\varphi(P_{\zeta_{0}})=\mathbb{E}\left[
(D\tilde{\varphi})(\zeta_{0})\cdot\eta\right]  +o(\left\vert \eta\right\vert
_{L^{2}\left(  \Omega\right)  }^{2})\text{, }\eta\in L^{2}\left(
\Omega,\mathcal{F},P;\mathbb{R}^{d}\right)  \text{.}%
\]
In Lions \cite{Lions} and Cardaliaguet \cite{Cardaliaguet}, it has been proved
that there exists a Borel function $h:\mathbb{R\rightarrow R}$, such that
$(D\tilde{\varphi})(\zeta_{0})=h(\zeta_{0})$ $P$-a.s. Note that $h(\zeta_{0})$
$P$-a.s. uniquely determined. Consequently, $h(y)$ is $P_{\zeta_{0}}(dy)$-a.e.
uniquely determined. We define%
\[
(\partial_{\mu}\varphi)(P_{\zeta_{0}},y):=h(y)\text{, }y\in\mathbb{R}\text{.}%
\]
Hence%

\[
\varphi(P_{\zeta_{0}+\eta})-\varphi(P_{\zeta_{0}})=\mathbb{E}\left[
(\partial_{\mu}\varphi)(P_{\zeta_{0}},\zeta_{0})\cdot\eta\right]
+o(\left\vert \eta\right\vert _{L^{2}\left(  \Omega\right)  }^{2})\text{,
}\left\vert \eta\right\vert _{L^{2}\left(  \Omega\right)  }^{2}\rightarrow
0\text{.}%
\]

\begin{example}
Given a function $\varphi(P_{\zeta})=g\left(  \mathbb{E}\left[  f(\zeta
)\right]  \right)  $, for $g,f\in C_{l,b}^{1}(\mathbb{R})$ and $\zeta\in
L^{2}\left(  \Omega,\mathcal{F},P;\mathbb{R}^{d}\right)  ,$ then
\begin{align*}
\mathbb{E}\left[  (\partial_{\mu}\varphi)(P_{\zeta},\zeta)\cdot\eta\right]
&  =\underset{\lambda\rightarrow0}{\lim}\left(  \varphi(P_{\zeta+\lambda\eta
})-\varphi(P_{\zeta})\right) \\
&  =\underset{\lambda\rightarrow0}{\lim}\frac{g\left(  \mathbb{E}\left[
f(\zeta+\lambda\eta)\right]  \right)  -g\left(  \mathbb{E}\left[
f(\zeta)\right]  \right)  }{\lambda}\\
&  =g^{\prime}\left(  \mathbb{E}\left[  f(\zeta)\right]  \right)
\mathbb{E}\left[  f^{\prime}(\zeta)\cdot\eta\right] \\
&  =\mathbb{E}\left[  g^{\prime}\mathbb{E}\left[  f(\zeta)\right]  f^{\prime
}(\zeta)\cdot\eta\right]  \text{, for all }\eta\in L^{2}\left(  \Omega
,\mathcal{F},P;\mathbb{R}^{d}\right)  \text{.}%
\end{align*}

\end{example}

Throughout this work, we will use also the following spaces:

\begin{itemize}
\item $S_{\mathbb{F}}^{2}([0,T])$ is the set of real valued $\mathbb{F}%
$-adapted continuous processes $(X(t))_{t\in\lbrack0,T]}$ such that
\[
{\Vert X\Vert}_{S_{\mathbb{F}}^{2}}:={\mathbb{E}}[\sup_{t\in\lbrack
0,T]}|X(t)|^{2}]<\infty.
\]

\item $L_{\mathbb{F}}^{2}([0,T])$ is the set of real valued $\mathbb{F}%
$-adapted processes $(Q(t))_{t\in\lbrack0,T]}$ such that
\[
\Vert Q\Vert_{L_{\mathbb{F}}^{2}}^{2}:={\mathbb{E}}[%
{\textstyle\int_{0}^{T}}
|Q(t)|^{2}dt]<\infty.
\]

\item $L^{2}(\mathcal{F}_{t})$ is the set of real valued square integrable
$\mathcal{F}_{t}$-measurable random variables.
\end{itemize}

\section{Solvability of the anticipated forward McKean-Vlasov equations}

Let us consider the following anticipated SDE for a given positive constant
$\delta$%

\begin{equation}
\left\{
\begin{array}
[c]{l}%
dX(t)=\sigma\left(  t,X(t),P_{X\left(  t+\delta\right)  }\right)
dB(t)+b\left(  t,X(t),P_{X\left(  t+\delta\right)  }\right)  dt\text{, }%
t\in\left[  0,T\right]  \text{,}\\
X\left(  0\right)  =x\in%
\mathbb{R}
^{d}\text{,}\\
X\left(  t\right)  =X\left(  T\right)  \text{, }t\geq T\text{.}%
\end{array}
\right.  \label{sde}%
\end{equation}
The functions $\sigma:\left[  0,T\right]  \times\Omega\times%
\mathbb{R}
^{d}\times\mathcal{P}_{2}\left(
\mathbb{R}
^{d}\right)  \rightarrow%
\mathbb{R}
^{d\times d}$ and $b:\left[  0,T\right]  \times\Omega\times%
\mathbb{R}
^{d}\times\mathcal{P}_{2}\left(
\mathbb{R}
^{d}\right)  \rightarrow%
\mathbb{R}
^{d}$ are progressively measurable and are assumed to satisfy the following
set of assumptions.\newline

\textbf{Assumptions (H.1)}

There exists $C>0$, such that

\begin{enumerate}
\item For all $t\in\left[  0,T\right]  $, $x,x^{\prime}\in%
\mathbb{R}
^{d},\mu,\mu^{\prime}\in\mathcal{P}_{2}\left(
\mathbb{R}
^{d}\right)  $%
\[
\left\vert \sigma\left(  t,x,\mu\right)  -\sigma\left(  t,x^{\prime}%
,\mu^{\prime}\right)  \right\vert +\left\vert b\left(  t,x,\mu\right)
-b\left(  t,x^{\prime},\mu^{\prime}\right)  \right\vert \leq C\left(
\left\vert x-x^{\prime}\right\vert +W_{2}\left(  \mu,\mu^{\prime}\right)
\right)  .
\]

\item For all $t\in\left[  0,T\right]  $, $x,x^{\prime}\in%
\mathbb{R}
^{d}$%
\[
\left\vert \sigma\left(  t,0,P_{0}\right)  \right\vert +\left\vert b\left(
t,0,P_{0}\right)  \right\vert \leq C\text{,}%
\]
where $P_{0}$ is the distribution law of zero, i.e., the Dirac measure with
mass at zero.
\end{enumerate}

\begin{remark}
Note that assumption (H.1) implies that the coefficients $b$ and $\sigma$ are
of linear growth. Indeed we have%
\begin{align*}
\left\vert \sigma\left(  t,x,\mu\right)  \right\vert  &  \leq\left\vert
\sigma\left(  t,0,P_{0}\right)  \right\vert +\left\vert \sigma\left(
t,x,\mu\right)  -\sigma\left(  t,0,P_{0}\right)  \right\vert \\
&  \leq C\left(  1+\left\vert x\right\vert +W_{2}\left(  \mu,P_{0}\right)
\right) \\
&  =C(1+\left\vert x\right\vert +(%
{\textstyle\int_{\mathbb{R}^{d}}}
\left\vert y\right\vert ^{2}\mu\left(  dy\right)  )^{\frac{1}{2}}),\left(
t,x,\mu\right)  \in\left[  0,T\right]  \times%
\mathbb{R}
^{d}\times\mathcal{P}_{2}\left(
\mathbb{R}
^{d}\right)  ,
\end{align*}
and a similar estimate holds for $b$.
\end{remark}

\begin{proposition}
\label{prop} Under the above assumption (H.1), there is some $\delta_{0}>0$,
such that for all $\delta\in\left(  0,\delta_{0}\right]  $, there exists a
unique solution $X\in S_{\mathbb{F}}^{2}\left(  \left[  0,T\right]  \right)  $
of SDE $\left(  \ref{sde}\right)  $; $\delta_{0}$ depends only on the
Lipschitz constant $C$ of the coefficients $b$ and $\sigma$ (see (H.1)) but
not on the coefficients themselves.
\end{proposition}

\noindent{Proof}\quad For $U\in S_{\mathbb{F}}^{2}\left(  \left[  0,T\right]
\right)  ,$ we can make the identification with the continuous process%
\[
\left(  U\left(  t\wedge T\right)  \right)  _{t\in\left[  0,T+\delta\right]
}\equiv(\left(  U\left(  t\right)  \right)  _{t\in\left[  0,T\right]
},U\left(  T\right)  )\in L_{\mathbb{F}}^{2}\left(  \left[  0,T\right]
\right)  \times L^{2}\left(  \mathcal{F}_{T}\right)  :=H.
\]
Given $U\in H$, we put%
\[
V\left(  t\right)  :=x+%
{\textstyle\int_{0}^{t}}
\sigma\left(  s,U\left(  s\right)  ,P_{U\left(  s+\delta\right)  }\right)
dB\left(  s\right)  +%
{\textstyle\int_{0}^{t}}
b\left(  s,U\left(  s\right)  ,P_{U\left(  s+\delta\right)  }\right)
ds,t\in\left[  0,T\right]  .
\]
Then $V\in S_{\mathbb{F}}^{2}\left(  \left[  0,T\right]  \right)  \subset H$
(with the above identification), and setting $\Phi\left(  U\right)  :=V$ we
define a mapping $\Phi:H\rightarrow H.$ Fixing $\beta>0$ ($\beta$ will be
specified later), we introduce the norm
\[
\left\Vert U\right\Vert _{-\beta}^{2}:=\mathbb{E[}e^{-\beta T}\left\vert
U\left(  T\right)  \right\vert ^{2}]+\tfrac{6}{7}\beta\mathbb{E[}%
{\textstyle\int_{0}^{T}}
e^{-\beta s}\left\vert U\left(  s\right)  \right\vert ^{2}ds]\text{, }U\in H.
\]
Obviously, $(H,\left\Vert \cdot\right\Vert _{-\beta})$ is a Banach space, and
the norm $\left\Vert \cdot\right\Vert _{-\beta}$ is equivalent to the norm
$\left\Vert \cdot\right\Vert _{0}$ (obtained from $\left\Vert \cdot\right\Vert
_{-\beta}$ by taking $\beta=0$)$.$ We are going to prove that $\Phi
:(H,\left\Vert \cdot\right\Vert _{-\beta})\rightarrow(H,\left\Vert
\cdot\right\Vert _{-\beta})$ is contracting. Indeed, we consider arbitrary
$U^{i}\in H$, $i=1,2,$ and we put $V^{i}:=\Phi\left(  U^{i}\right)  $,
$i=1,2$. Let $\bar{U}:=U^{1}-U^{2}$ and $\bar{V}:=V^{1}-V^{2}.$ Then, applying
Itô's formula to $(e^{-\beta t}\left\vert \bar{V}\left(  t\right)  \right\vert
^{2})_{t\geq0}$, we get from the assumptions (H.1)
\[%
\begin{array}
[c]{l}%
\mathbb{E[}e^{-\beta t}\left\vert \bar{V}\left(  t\right)  \right\vert
^{2}]+\mathbb{E[}%
{\textstyle\int_{0}^{t}}
\beta e^{-\beta s}\left\vert \bar{V}\left(  s\right)  \right\vert ^{2}ds]\\
=2\mathbb{E[}%
{\textstyle\int_{0}^{t}}
e^{-\beta s}\bar{V}\left(  s\right)  (b(s,U^{1}\left(  s\right)
,P_{U^{1}\left(  s+\delta\right)  })-b(s,U^{2}\left(  s\right)  ,P_{U^{2}%
\left(  s+\delta\right)  }))ds]\\
+\mathbb{E[}%
{\textstyle\int_{0}^{t}}
e^{-\beta s}|\sigma(s,U^{1}\left(  s\right)  ,P_{U^{1}\left(  s+\delta\right)
})-\sigma(s,U^{2}\left(  s\right)  ,P_{U^{2}\left(  s+\delta\right)  }%
)|^{2}ds]\\
\leq C\mathbb{E[}%
{\textstyle\int_{0}^{t}}
e^{-\beta s}\left\vert \bar{V}\left(  s\right)  \right\vert (\left\vert
\bar{U}\left(  s\right)  \right\vert +W_{2}(P_{U^{1}\left(  s+\delta\right)
},P_{U^{2}\left(  s+\delta\right)  }))ds]\\
+C\mathbb{E[}%
{\textstyle\int_{0}^{t}}
e^{-\beta s}(\left\vert \bar{U}\left(  s\right)  \right\vert +W_{2}%
(P_{U^{1}\left(  s+\delta\right)  },P_{U^{2}\left(  s+\delta\right)  }%
))^{2}ds]\\
\leq C\mathbb{E[}%
{\textstyle\int_{0}^{t}}
e^{-\beta s}\left\vert \bar{V}\left(  s\right)  \right\vert ^{2}%
ds]+C\mathbb{E[}%
{\textstyle\int_{0}^{t}}
e^{-\beta s}\left\vert \bar{U}\left(  s\right)  \right\vert ^{2}ds]\\
+C\mathbb{E[}%
{\textstyle\int_{0}^{t}}
e^{-\beta s}\left\vert \bar{U}\left(  s+\delta\right)  \right\vert
^{2}ds],t\in\left[  0,T\right]  .
\end{array}
\]
Indeed, we recall that%
\[
W_{2}^{2}(P_{U^{1}\left(  s+\delta\right)  },P_{U^{2}\left(  s+\delta\right)
})\leq\mathbb{E[}\left\vert U^{1}\left(  s+\delta\right)  -U^{2}\left(
s+\delta\right)  \right\vert ^{2}]=\mathbb{E[}\left\vert \bar{U}\left(
s+\delta\right)  \right\vert ^{2}].
\]
Hence for $t=T$, we have%
\[%
\begin{array}
[c]{l}%
\mathbb{E[}e^{-\beta T}\left\vert \bar{V}\left(  T\right)  \right\vert
^{2}]+\mathbb{E[}%
{\textstyle\int_{0}^{T}}
\beta e^{-\beta s}\left\vert \bar{V}\left(  s\right)  \right\vert ^{2}ds]\\
\leq C\mathbb{E[}%
{\textstyle\int_{0}^{T}}
e^{-\beta s}\left\vert \bar{V}\left(  s\right)  \right\vert ^{2}%
ds]+C\mathbb{E[}%
{\textstyle\int_{0}^{T}}
e^{-\beta s}\left\vert \bar{U}\left(  s\right)  \right\vert ^{2}ds]\\
+Ce^{\beta\delta}\mathbb{E[}%
{\textstyle\int_{0}^{T}}
e^{-\beta s}\left\vert \bar{U}\left(  s\right)  \right\vert ^{2}%
ds]+Ce^{\beta\delta}\delta\mathbb{E[}e^{-\beta T}\left\vert \bar{U}\left(
T\right)  \right\vert ^{2}].
\end{array}
\]
We seek suitable $\beta>0$, $\delta>0$ with $\delta\leq\frac{1}{\beta}$, i.e.,
$\beta\delta\leq1$, in order to estimate%
\[%
\begin{array}
[c]{l}%
\mathbb{E[}e^{-\beta T}\left\vert \bar{V}\left(  T\right)  \right\vert
^{2}]+\beta\mathbb{E[}%
{\textstyle\int_{0}^{T}}
e^{-\beta s}\left\vert \bar{V}\left(  s\right)  \right\vert ^{2}ds]\\
\leq Ce\delta\mathbb{E[}e^{-\beta T}\left\vert \bar{U}\left(  T\right)
\right\vert ^{2}]+C\mathbb{E[}%
{\textstyle\int_{0}^{T}}
e^{-\beta s}\left\vert \bar{V}\left(  s\right)  \right\vert ^{2}ds]\\
+C\left(  1+e\right)  \mathbb{E[}%
{\textstyle\int_{0}^{T}}
e^{-\beta s}\left\vert \bar{U}\left(  s\right)  \right\vert ^{2}ds].
\end{array}
\]
Choosing $\beta:=7C$, $\delta_{0}:=\tfrac{1}{7C}(=\frac{1}{\beta})$, we have
for all $\delta\in\left(  0,\delta_{0}\right)  $:%
\[%
\begin{array}
[c]{l}%
\mathbb{E[}e^{-\beta T}\left\vert \bar{V}\left(  T\right)  \right\vert
^{2}]+6C\mathbb{E[}%
{\textstyle\int_{0}^{T}}
e^{-\beta s}\left\vert \bar{V}\left(  s\right)  \right\vert ^{2}ds]\\
\leq\tfrac{2}{3}(\mathbb{E[}e^{-\beta T}\left\vert \bar{U}\left(  T\right)
\right\vert ^{2}]+6C\mathbb{E[}%
{\textstyle\int_{0}^{T}}
e^{-\beta s}\left\vert \bar{U}\left(  s\right)  \right\vert ^{2}ds]).
\end{array}
\]
Then%
\[
\left\Vert \bar{V}\right\Vert _{-\beta}\leq(\tfrac{2}{3})^{\frac{1}{2}%
}\left\Vert \bar{U}\right\Vert _{-\beta},
\]
i.e.,%
\[
\left\Vert \Phi\left(  U^{1}\right)  -\Phi\left(  U^{2}\right)  \right\Vert
_{-\beta}\leq(\tfrac{2}{3})^{\frac{1}{2}}\left\Vert U^{1}-U^{2}\right\Vert
_{-\beta},\text{ for all }U^{1},U^{2}\in H.
\]
This proves that $\Phi:(H,\left\Vert \cdot\right\Vert _{-\beta})\rightarrow
(H,\left\Vert \cdot\right\Vert _{-\beta})$ is a contraction on the Banach
space $\left(  H,\left\Vert \cdot\right\Vert _{-\beta}\right)  $. Hence, there
is a unique fixed point $X\in H,$ such that $X=\Phi\left(  X\right)  ,$ i.e.,%
\[
X\left(  t\right)  =x+%
{\textstyle\int_{0}^{t}}
\sigma\left(  s,X\left(  s\right)  ,P_{X\left(  s+\delta\right)  }\right)
dB\left(  s\right)  +%
{\textstyle\int_{0}^{t}}
b\left(  s,X\left(  s\right)  ,P_{X\left(  s+\delta\right)  }\right)  ds,
\]
$v\left(  dt\right)  $-a.e. on $\left[  0,T\right]  $, $P$-a.s., with
$v\left(  dt\right)  =I_{\left[  0,T\right]  }\left(  t\right)  dt+P_{T}%
\left(  dt\right)  $ (Recall the definition of $H$). For a $v\otimes
P-$modification of $X,$ also denoted by $X$, we have $X\in S_{\mathbb{F}}%
^{2}\left(  \left[  0,T\right]  \right)  $ and
\[
X\left(  t\right)  =x+%
{\textstyle\int_{0}^{t}}
\sigma\left(  s,X\left(  s\right)  ,P_{X\left(  s+\delta\right)  }\right)
dB\left(  s\right)  +%
{\textstyle\int_{0}^{t}}
b\left(  s,X\left(  s\right)  ,P_{X\left(  s+\delta\right)  }\right)
ds,t\in\left[  0,T\right]  \text{ }P\text{-a.s.}%
\]
{ }$\square$

\section{Pontryagin's stochastic maximum principle}

Let us introduce now our stochastic control problem.

\subsection{Controlled stochastic differential equation}

As control state space we consider a bounded convex subset $U$ of $%
\mathbb{R}
^{d}$. A process $u=\left(  u(t)\right)  _{t\in\left[  0,T\right]  }:\left[
0,T\right]  \times\Omega\rightarrow U$ which is progressively measurable is
called an admissible control; $\mathcal{U=}L_{\mathbb{F}}^{0}(\left[
0,T\right]  ;U)$ is the set of all admissible controls. The dynamics of our
controlled system are driven by functions $\sigma:\left[  0,T\right]
\times\Omega\times%
\mathbb{R}
^{d}\times\mathcal{P}_{2}\left(
\mathbb{R}
^{d}\right)  \times U\rightarrow%
\mathbb{R}
^{d\times d}$, $b:\left[  0,T\right]  \times\Omega\times%
\mathbb{R}
^{d}\times\mathcal{P}_{2}\left(
\mathbb{R}
^{d}\right)  \times U\rightarrow%
\mathbb{R}
^{d}$.\newline

\textbf{Assumptions (H.2)}

The coefficients $\sigma$ and $b$ are supposed to be continuous on $\left[
0,T\right]  \times\Omega\times%
\mathbb{R}
^{d}\times\mathcal{P}_{2}\left(
\mathbb{R}
^{d}\right)  \times U$ and Lipschitz on $%
\mathbb{R}
^{d}\times\mathcal{P}_{2}\left(
\mathbb{R}
^{d}\right)  ,$ uniformly w.r.t. $u\in U$ and $\omega\in\Omega$ i.e., there is
some $C>0$, such that for all $\left(  x,\mu\right)  ,\left(  x^{\prime}%
,\mu^{\prime}\right)  \in%
\mathbb{R}
^{d}\times\mathcal{P}_{2}\left(
\mathbb{R}
^{d}\right)  $, $u\in U$, we have%

\begin{align*}
\left\vert \sigma\left(  t,x,\mu,u\right)  -\sigma\left(  t,x^{\prime}%
,\mu^{\prime},u\right)  \right\vert  &  \leq C\left(  \left\vert x-x^{\prime
}\right\vert +W_{2}\left(  \mu^{\prime},\mu\right)  \right)  ,\\
\left\vert b\left(  t,x,\mu,u\right)  -b\left(  t,x^{\prime},\mu^{\prime
},u\right)  \right\vert  &  \leq C\left(  \left\vert x-x^{\prime}\right\vert
+W_{2}\left(  \mu^{\prime},\mu\right)  \right)  .
\end{align*}
On the other hand, from the continuity of the coefficients on $\left[
0,T\right]  \times\Omega\times%
\mathbb{R}
^{d}\times\mathcal{P}_{2}\left(
\mathbb{R}
^{d}\right)  \times U$, we have
\[
\left\vert \sigma\left(  t,0,P_{0},u\right)  \right\vert +\left\vert b\left(
t,0,P_{0},u\right)  \right\vert \leq C\text{, for all }u\in U.
\]
This shows that, for every $u\in\mathcal{U}:=L_{\mathbb{F}}^{0}\left(  \left[
0,T\right]  ;U\right)  $ and $\omega\in\Omega$, the coefficients $\sigma$ and
$b$ satisfy the assumptions (H.1). Thus, for $\delta_{0}>0$ from Proposition
\ref{prop}, for all $u\in\mathcal{U}$; $x\in%
\mathbb{R}
^{d}$, there is a unique solution $X^{u}\left(  t\right)  \in S_{\mathbb{F}%
}^{2}\left(  \left[  0,T\right]  ;%
\mathbb{R}
^{d}\right)  $ of the equation%

\[
X^{u}\left(  t\right)  =x+%
{\textstyle\int_{0}^{t}}
\sigma\left(  s,X^{u}\left(  s\right)  ,P_{X^{u}(s+\delta)},u(s)\right)
dB\left(  s\right)  +%
{\textstyle\int_{0}^{t}}
b\left(  s,X^{u}\left(  s\right)  ,P_{X^{u}(s+\delta)},u\left(  s\right)
\right)  ds,t\in\left[  0,T\right]  .
\]

\subsection{Cost functional}

Let us endow our control problem with a terminal cost $g:\left[  0,T\right]
\times\Omega\times%
\mathbb{R}
^{d}\times\mathcal{P}_{2}\left(
\mathbb{R}
^{d}\right)  \rightarrow%
\mathbb{R}
$, and a running cost $l:\left[  0,T\right]  \times\Omega\times%
\mathbb{R}
^{d}\times\mathcal{P}_{2}\left(
\mathbb{R}
^{d}\right)  \times U\rightarrow%
\mathbb{R}
$.\newline

\textbf{Assumptions (H.3)}

We suppose that $g:\left[  0,T\right]  \times\Omega\times%
\mathbb{R}
^{d}\times\mathcal{P}_{2}\left(
\mathbb{R}
^{d}\right)  \rightarrow%
\mathbb{R}
$ is continuous and satisfies a linear growth assumption: For some constant
$C>0$,%

\[
\left\vert g\left(  x,\mu\right)  \right\vert \leq C(1+\left\vert x\right\vert
+(%
{\textstyle\int_{\mathbb{R}^{d}}}
\left\vert y\right\vert ^{2}\mu\left(  dy\right)  )^{\frac{1}{2}})\text{,
}\left(  x,\mu\right)  \in%
\mathbb{R}
^{d}\times\mathcal{P}_{2}\left(
\mathbb{R}
^{d}\right)  \text{.}%
\]
Let $l:\left[  0,T\right]  \times\Omega\times%
\mathbb{R}
^{d}\times\mathcal{P}_{2}\left(
\mathbb{R}
^{d}\right)  \times U\rightarrow%
\mathbb{R}
$ be continuous and such that, for some $C>0$, for all $\left(  x,\mu
,u\right)  \in%
\mathbb{R}
^{d}\times\mathcal{P}_{2}\left(
\mathbb{R}
^{d}\right)  \times U$, we have%

\[
\left\vert l\left(  x,\mu,u\right)  \right\vert \leq C(1+\left\vert
x\right\vert +(%
{\textstyle\int_{\mathbb{R}^{d}}}
\left\vert y\right\vert ^{2}\mu\left(  dy\right)  )^{\frac{1}{2}})\text{.}%
\]
For any admissible control $u$, we define the performance functional:%

\[
J\left(  u\right)  :=\mathbb{E[}g\left(  X^{u}\left(  T\right)  ,P_{X^{u}%
(T)}\right)  +%
{\textstyle\int_{0}^{T}}
l\left(  s,X^{u}\left(  s\right)  ,P_{X^{u}(s+\delta)},u\left(  s\right)
\right)  ds]\text{.}%
\]
A control process $u^{\ast}\in\mathcal{U}$ is called optimal, if%

\[
J\left(  u^{\ast}\right)  \leq J\left(  u\right)  \text{, for all }%
u\in\mathcal{U}\text{.}%
\]
Let us suppose that there is an optimal control $u^{\ast}\in\mathcal{U}$. Our
objective is to characterise the optimal control. For this let us assume some
additional assumptions.\newline

\textbf{Assumptions (H.4)}

Let $U$ be convex (and, hence, $\mathcal{U}$ is convex). The functions
$\sigma\left(  \cdot,\cdot,\cdot,u\right)  $, $b\left(  \cdot,\cdot
,\cdot,u\right)  $, $l\left(  \cdot,\cdot,\cdot,u\right)  $ and $g\left(
\cdot,\cdot\right)  $ are continuously differentiable over $%
\mathbb{R}
^{d}\times\mathcal{P}_{2}\left(
\mathbb{R}
^{d}\right)  \times U$ with bounded derivatives.\newline

Given an arbitrary but fixed control $u\in\mathcal{U}$, we define
\[%
\begin{array}
[c]{ll}%
u^{\theta}:=u^{\ast}+\theta\left(  u-u^{\ast}\right)  , & \theta\in\left[
0,1\right]  .
\end{array}
\]
Note that, thanks to the convexity of $U$ and $\mathcal{U}$, also $u^{\theta
}\in\mathcal{U},\theta\in\left[  0,1\right]  $. We denote by $X^{\theta
}:=X^{u^{\theta}}$ and by $X^{\ast}:=X^{u^{\ast}}$ the solution processes
corresponding to $u^{\theta}$ and $u^{\ast},$\ respectively. For simplicity of
the computations, we set $d=1.$

\subsection{Variational SDE}

Given $u^{\ast}\in\mathcal{U}$ and the associated controlled state process
$X^{\ast}$, let $Y=\left(  Y(t)\right)  _{t\in\left[  0,T\right]  }\in
S_{\mathbb{F}}^{2}\left(  \left[  0,T\right]  \right)  \left(  :=S_{\mathbb{F}%
}^{2}\left(  \left[  0,T\right]  ;%
\mathbb{R}
\right)  \right)  $ be the unique solution of the following SDE%

\begin{equation}
\left\{
\begin{array}
[c]{ll}%
Y\left(  t\right)  & =%
{\textstyle\int_{0}^{t}}
\{(\partial_{x}\sigma)(s,X^{\ast}\left(  s\right)  ,P_{X^{\ast}(s+\delta
)},u^{\ast}\left(  s\right)  )Y\left(  s\right) \\
& +\mathbb{\tilde{E}[(}\partial_{\mu}\sigma)(s,X^{\ast}\left(  s\right)
,P_{X^{\ast}(s+\delta)},\tilde{X}^{\ast}(s+\delta))\tilde{Y}(s+\delta)]\\
& +(\partial_{u}\sigma)\left(  s,X^{\ast}\left(  s\right)  ,P_{X^{\ast
}(s+\delta)},u^{\ast}\left(  s\right)  \right)  (u\left(  s\right)  -u^{\ast
}\left(  s\right)  )\}dB\left(  s\right) \\
& +%
{\textstyle\int_{0}^{t}}
\{(\partial_{x}b)(s,X^{\ast}\left(  s\right)  ,P_{X^{\ast}(s+\delta)},u^{\ast
}\left(  s\right)  )Y(s)\\
& +\mathbb{\tilde{E}[(}\partial_{\mu}b)(s,X^{\ast}\left(  s\right)
,P_{X^{\ast}(s+\delta)},\tilde{X}^{\ast}(s+\delta))\tilde{Y}(s+\delta)]\\
& +(\partial_{u}b)\left(  s,X^{\ast}\left(  s\right)  ,P_{X^{\ast}(s+\delta
)},u^{\ast}\left(  s\right)  \right)  (u\left(  s\right)  -u^{\ast}\left(
s\right)  )\}ds,t\in\lbrack0,T],\\
Y\left(  0\right)  & =0\text{, }Y(t)=Y(T)\text{, }t\geq T\text{.}%
\end{array}
\right.  \label{y}%
\end{equation}

\begin{remark}
Note that SDE $\left(  \ref{y}\right)  $ is obtained by formal differentiation
of equation $\left(  \ref{x}\right)  $ (with $u=u^{\theta}$) at $\theta=0$.
\end{remark}

From the previous section, we have the existence and the uniqueness of a
solution for all $\delta\in\left(  0,\delta_{0}^{\prime}\right]  ;0<$
$\delta_{0}^{\prime}\leq\delta_{0}$ small enough.\newline

Indeed, equation $\left(  \ref{y}\right)  $ is of the form%
\[%
\begin{array}
[c]{ll}%
Y\left(  t\right)  & =%
{\textstyle\int_{0}^{t}}
\alpha_{1}(s)Y\left(  s\right)  ds+%
{\textstyle\int_{0}^{t}}
\alpha_{2}(s)ds+%
{\textstyle\int_{0}^{t}}
\mathbb{\tilde{E}[}\beta_{1}(s)\tilde{Y}(s+\delta)]ds\\
& +%
{\textstyle\int_{0}^{t}}
\alpha_{3}(s)Y(s)dB(s)+%
{\textstyle\int_{0}^{t}}
\alpha_{4}(s)dB(s)+%
{\textstyle\int_{0}^{t}}
\mathbb{\tilde{E}[}\beta_{2}(s)\tilde{Y}(s+\delta)]dB(s),t\in\left[
0,T\right]  ,
\end{array}
\]
where $\alpha_{i}:\left[  0,T\right]  \times\Omega\rightarrow%
\mathbb{R}
,i=1,2,3,4,$ are bounded progressively measurable processes and $\beta
_{j}:\left[  0,T\right]  \times\Omega\rightarrow%
\mathbb{R}
,j=1,2,$ two bounded $(\mathcal{F}_{t}\otimes\mathcal{\tilde{F}}%
_{T}\mathcal{)}$- progressively measurable processes. With the method used in
the proof of Proposition \ref{prop}, we get the existence of $\delta
_{0}^{^{\prime}}\in\left(  0,\delta_{0}\right]  $ stated above. It turns out
that $Y(t)$ is the $L^{2}$-derivative of $X^{\theta}(t)$ w.r.t. $\theta$ at
$\theta=0.$ More precisely, the following property holds.

\begin{lemma}
\label{lem} $\mathbb{E[}\sup_{t\in\left[  0,T\right]  }|Y(t)-\tfrac{X^{\theta
}(t)-X^{\ast}(t)}{\theta}|^{2}]\rightarrow0$ as $\left(  \theta\rightarrow
0\right)  $.
\end{lemma}

\noindent{Proof}\quad The proof is obtained with standard computations. For
the sake of completeness, we give details in the Appendix.{ }$\square$

\subsection{Variational inequality}

We know that if $u^{\ast}$ is an optimal control, we have $J(u^{\ast})\leq
J(u^{\theta})$, for all $\theta\in\left[  0,1\right]  $, i.e.,%

\begin{equation}
0\leq\underset{\theta\rightarrow0}{\underline{\lim}}\tfrac{J(u^{\theta
})-J(u^{\ast})}{\theta}\text{.} \label{j}%
\end{equation}

\begin{lemma}
Under assumptions (H.3), (H.4), Lemma \ref{lem} and inequality $\left(
\ref{j}\right)  $, we have%
\begin{equation}%
\begin{array}
[c]{l}%
0\leq\mathbb{E[\{(}\partial_{x}g)(X^{\ast}(T),P_{X^{\ast}(T)})+\mathbb{\tilde
{E}[(}\partial_{\mu}g)(\tilde{X}^{\ast}(T),P_{X^{\ast}(T)},X^{\ast}(T))]\\
+%
{\textstyle\int_{T-\delta}^{T}}
\mathbb{\tilde{E}[(}\partial_{\mu}l)(t,\tilde{X}^{\ast}(t),P_{X^{\ast}%
(T)},X^{\ast}(T),\tilde{u}^{\ast}(t))dt]\}Y(T)]\\
+\mathbb{E[}%
{\textstyle\int_{0}^{T}}
\{(\partial_{x}l)(t,X^{\ast}(t),P_{X^{\ast}(t+\delta)},u^{\ast}(t))\\
+\mathbb{\tilde{E}[(}\partial_{\mu}l)(t-\delta,\tilde{X}^{\ast}(t-\delta
),P_{X^{\ast}(t)},x^{\ast}(t),\tilde{u}^{\ast}(t-\delta))]I_{\left[
\delta,T\right]  }(t)\}Y(t)dt]\\
+\mathbb{E[}%
{\textstyle\int_{0}^{T}}
(\partial_{u}l)(t,X^{\ast}(t),P_{X^{\ast}(t+\delta)},u^{\ast}(t))(u(t)-u^{\ast
}(t))dt].
\end{array}
\label{B}%
\end{equation}

\end{lemma}

\noindent{Proof}\quad From the definition of $J(u^{\ast})$, we have%
\[%
\begin{array}
[c]{ll}%
\underset{\theta\rightarrow0}{\underline{\lim}}\tfrac{1}{\theta}(J(u^{\theta
})-J(u^{\ast})) & =\mathbb{E[}(\partial_{x}g)\left(  X^{\ast}(T),P_{X^{\ast
}(T)})Y(T)\right)  \\
& +\mathbb{\tilde{E}[(}\partial_{\mu}g)(X^{\ast}(T),P_{X^{\ast}(T)},\tilde
{X}^{\ast}(T))\tilde{Y}(T)]\\
& +%
{\textstyle\int_{0}^{T}}
\{(\partial_{x}l)(t,X^{\ast}(t),P_{X^{\ast}(t+\delta)},u^{\ast}(t))Y(t)\\
& +\mathbb{\tilde{E}[}(\partial_{\mu}l)(t,X^{\ast}(t),P_{X^{\ast}(t+\delta
)},\tilde{X}^{\ast}(t+\delta),u^{\ast}(t))\tilde{Y}(t+\delta)]\\
& +(\partial_{u}l)(t,X^{\ast}(t),P_{X^{\ast}(t+\delta)},u^{\ast}%
(t))(u(t)-u^{\ast}(t))\}dt]\text{,}%
\end{array}
\]
the fact that
\[
\mathbb{E[}\underset{t\in\left[  0,T\right]  }{\sup}|Y(t)-\tfrac{X^{\theta
}(t)-X^{\ast}(t)}{\theta}|^{2}]\rightarrow0\text{ as }\theta\rightarrow
0\text{,}%
\]
\bigskip and by repeating previous arguments, $(\ref{B})$ is obtained.$\qquad
\square$

\subsection{Adjoint processes}

Let us first recall the equation satisfied by the derivative process
\[
\left\{
\begin{array}
[c]{ll}%
dY\left(  t\right)  & =\{\left(  \partial_{x}\sigma\right)  (t)Y\left(
t\right)  +\mathbb{\tilde{E}[}\left(  \partial_{\mu}\sigma\right)
(t)\tilde{Y}\left(  t+\delta\right)  ]\\
& +\left(  \partial_{u}\sigma\right)  (t)\left(  u\left(  t\right)  -u^{\ast
}\left(  t\right)  \right)  \}dB\left(  t\right) \\
& +\{\left(  \partial_{x}b\right)  (t)Y\left(  t\right)  +\mathbb{\tilde{E}%
[}\left(  \partial_{\mu}b\right)  (t)\tilde{Y}\left(  t+\delta\right)  ]\\
& +\left(  \partial_{u}b\right)  (t)\left(  u\left(  t\right)  -u^{\ast
}\left(  t\right)  \right)  \}dt,\\
Y\left(  0\right)  & =0,
\end{array}
\right.
\]
where for notational convenient, we have used the short hand notations%
\[%
\begin{array}
[c]{l}%
\left(  \partial_{x}\sigma\right)  (t,X^{^{\ast}}\left(  t\right)
,P_{X^{^{\ast}}\left(  t+\delta\right)  },u^{\ast}\left(  t\right)  )=:\left(
\partial_{x}\sigma\right)  (t),\\
\mathbb{\tilde{E}[}\left(  \partial_{\mu}\sigma\right)  (t,X^{^{\ast}}\left(
t\right)  ,P_{X^{^{\ast}}\left(  t+\delta\right)  },\tilde{X}^{\ast}\left(
t+\delta\right)  ,u^{\ast}(t))\tilde{Y}\left(  t+\delta\right)
]=:\mathbb{\tilde{E}[}\left(  \partial_{\mu}\sigma\right)  (t)\tilde{Y}\left(
t+\delta\right)  ],
\end{array}
\]
and similarly.In order to determine the adjoint backward equation, we suppose
that it has the form%

\begin{equation}
\left\{
\begin{array}
[c]{ll}%
dp(t) & =-\alpha(t)dt+q(t)dB(t),t\in\left[  0,T\right]  ,\\
p(T), &
\end{array}
\right.  \label{bac}%
\end{equation}
for some adapted process $\alpha$ and terminal value $p(T)$ which we have to
determine. Applying Itô's formula to $p\left(  t\right)  Y\left(  t\right)  ,$
we obtain%

\begin{equation}%
\begin{array}
[c]{l}%
d\mathbb{E}\left[  p\left(  t\right)  Y\left(  t\right)  \right]
=\{\mathbb{E[}\left(  \partial_{x}b\right)  (t)p\left(  t\right)  Y\left(
t\right)  ]\\
+\mathbb{E[\tilde{E}[}\left(  \partial_{\mu}b\right)  (t)\tilde{Y}\left(
t+\delta\right)  ]p\left(  t\right)  ]+\mathbb{E}\left[  \left(  \partial
_{u}b\right)  \left(  t\right)  \left(  u\left(  t\right)  -u^{\ast}\left(
t\right)  \right)  p\left(  t\right)  \right]  \}dt\\
-\mathbb{E}\left[  \alpha(t)Y(t)\right]  dt+\{\mathbb{E[}\left(  \partial
_{x}\sigma\right)  (t)q(t)Y\left(  t\right)  ]+\mathbb{E[\tilde{E}[}\left(
\partial_{\mu}\sigma\right)  (t)\tilde{Y}(t+\delta)]q(t)]\\
+\mathbb{E[}(\partial_{u}\sigma)(t)\left(  u\left(  t\right)  -u^{\ast}\left(
t\right)  \right)  q(t)]\}dt,
\end{array}
\label{7}%
\end{equation}
with $\mathbb{E}\left[  Y\left(  0\right)  p\left(  0\right)  \right]  =0.$

We have that%
\begin{equation}%
\begin{array}
[c]{ll}%
\mathbb{E[\tilde{E}[}f(\xi,\tilde{\eta})]] & =%
{\textstyle\int_{\Omega}}
(%
{\textstyle\int_{\tilde{\Omega}}}
f(\xi(\omega),\tilde{\eta}(\tilde{\omega}))\tilde{P}(d\tilde{\omega
}))P(d\omega)\\
& =%
{\textstyle\int_{\Omega}}
(%
{\textstyle\int_{\mathbb{R}}}
f(\xi,y)\tilde{P}_{\tilde{\eta}}(dy))dP\\
& =%
{\textstyle\int_{\Omega}}
(%
{\textstyle\int_{\mathbb{R}}}
f(x,y)\tilde{P}_{\tilde{\eta}}(dy))P_{\xi}(dx)\\
& =%
{\textstyle\int_{\mathbb{R}}}
(%
{\textstyle\int_{\Omega}}
f(x,y)dP)P_{\xi}(dx)\\
& =%
{\textstyle\int_{\tilde{\Omega}}}
(%
{\textstyle\int_{\Omega}}
f(\tilde{\xi},\eta)dP)d\tilde{P}\\
& =\mathbb{\tilde{E}[E[}f(\tilde{\xi},\eta)]].
\end{array}
\label{eetilde}%
\end{equation}
Using the above computations, we obtain%
\begin{equation}%
\begin{array}
[c]{c}%
\mathbb{E[\tilde{E}[}\left(  \partial_{\mu}b\right)  (t,X^{\ast}\left(
t\right)  ,P_{X^{^{\ast}}\left(  t+\delta\right)  },\tilde{X}^{\ast}\left(
t+\delta\right)  ,u^{\ast}\left(  t\right)  )\tilde{Y}\left(  t+\delta\right)
]p\left(  t\right)  ]\\
=\mathbb{E[\tilde{E}[}\left(  \partial_{\mu}b\right)  (t,\tilde{X}^{\ast
}\left(  t\right)  ,P_{X^{^{\ast}}\left(  t+\delta\right)  },X^{\ast}\left(
t+\delta\right)  ,\tilde{u}^{\ast}\left(  t\right)  )\tilde{p}\left(
t\right)  ]Y\left(  t+\delta\right)  ],
\end{array}
\label{8}%
\end{equation}
and similarly%
\begin{equation}%
\begin{array}
[c]{c}%
\mathbb{E[\tilde{E}[}\left(  \partial_{\mu}\sigma\right)  (t,X^{\ast}\left(
t\right)  ,P_{X^{^{\ast}}\left(  t+\delta\right)  },\tilde{X}^{^{\ast}}\left(
t+\delta\right)  ,u^{\ast}\left(  t\right)  )\tilde{Y}(t+\delta)]q(t)]\\
=\mathbb{E[\tilde{E}[}\left(  \partial_{\mu}\sigma\right)  (t,\tilde{X}^{\ast
}\left(  t\right)  ,P_{X^{^{\ast}}\left(  t+\delta\right)  },X^{^{\ast}%
}\left(  t+\delta\right)  ,\tilde{u}^{\ast}\left(  t\right)  )\tilde
{q}(t)]Y(t+\delta)]\text{,}%
\end{array}
\label{9}%
\end{equation}
Substituting $\left(  \ref{8}\right)  -\left(  \ref{9}\right)  $ into $\left(
\ref{7}\right)  ,$ we get%
\begin{equation}%
\begin{array}
[c]{ll}%
\mathbb{E[}p\left(  T\right)  Y\left(  T\right)  ] & =%
{\textstyle\int_{0}^{T}}
\mathbb{E[}((\partial_{x}b)(t)p(t)+(\partial_{x}\sigma)(t)q(t))Y(t)]dt\\
& +%
{\textstyle\int_{0}^{T}}
\mathbb{E[\tilde{E}[(}\left(  \partial_{\mu}b\right)  (t)\tilde{p}\left(
t\right)  +\left(  \partial_{\mu}\sigma\right)  (t)\tilde{q}\left(  t\right)
)]Y\left(  t+\delta\right)  ]dt-%
{\textstyle\int_{0}^{T}}
\mathbb{E[}\alpha(t)Y(t)]dt\\
& +%
{\textstyle\int_{0}^{T}}
\mathbb{E[(}(\partial_{u}b)(t)p(t)+(\partial_{u}\sigma)(t)q(t))\left(
u\left(  t\right)  -u^{\ast}\left(  t\right)  \right)  ]dt.
\end{array}
\label{10}%
\end{equation}
As $X^{\ast}\left(  t\right)  =X^{\ast}\left(  T\right)  ,$ $Y(t)=Y(T),t\geq
T$, we get
\begin{equation}%
\begin{array}
[c]{l}%
{\textstyle\int_{0}^{T}}
\mathbb{E[\tilde{E}[}\left(  \partial_{\mu}b\right)  (t)\tilde{p}\left(
t\right)  ]Y\left(  t+\delta\right)  ]dt\\
=\mathbb{E[(}%
{\textstyle\int_{T-\delta}^{T}}
\mathbb{\tilde{E}[}\left(  \partial_{\mu}b\right)  (t)\tilde{p}\left(
t\right)  ]dt)Y\left(  T\right)  ]\\
+\mathbb{E[}%
{\textstyle\int_{0}^{T}}
\mathbb{\tilde{E}[}\left(  \partial_{\mu}b\right)  (t-\delta)\tilde{p}\left(
t-\delta\right)  ]I_{\left[  \delta,T\right]  }\left(  t\right)  Y\left(
t\right)  dt].
\end{array}
\label{11}%
\end{equation}
Analogously,%
\begin{equation}%
\begin{array}
[c]{l}%
{\textstyle\int_{0}^{T}}
\mathbb{E[\tilde{E}[}\left(  \partial_{\mu}\sigma\right)  (t)\tilde{q}\left(
t\right)  ]Y\left(  t+\delta\right)  ]dt\\
=\mathbb{E[(}%
{\textstyle\int_{T-\delta}^{T}}
\mathbb{\tilde{E}[}\left(  \partial_{\mu}\sigma\right)  (t)\tilde{q}\left(
t\right)  ]dt)Y\left(  T\right)  ]\\
+\mathbb{E[}%
{\textstyle\int_{0}^{T}}
\mathbb{\tilde{E}[}\left(  \partial_{\mu}b\right)  (t-\delta)\tilde{q}\left(
t-\delta\right)  ]I_{\left[  \delta,T\right]  }\left(  t\right)  Y\left(
t\right)  dt].
\end{array}
\label{12}%
\end{equation}
Combining $\left(  \ref{10}\right)  -\left(  \ref{12}\right)  $, we obtain%
\begin{equation}%
\begin{array}
[c]{l}%
\mathbb{E[(}p(T)-%
{\textstyle\int_{T-\delta}^{T}}
\mathbb{\tilde{E}[}\left(  \partial_{\mu}b\right)  (t)\tilde{p}\left(
t\right)  +\left(  \partial_{\mu}\sigma\right)  (t)\tilde{q}\left(  t\right)
]dt)Y(T)]\\
=\mathbb{E[}%
{\textstyle\int_{0}^{T}}
\{(\partial_{x}b)(t)p(t)+(\partial_{x}\sigma)(t)q(t)\\
+\mathbb{\tilde{E}[}\left(  \partial_{\mu}b\right)  (t-\delta)\tilde{p}\left(
t-\delta\right)  +\left(  \partial_{\mu}\sigma\right)  (t-\delta)\tilde
{q}\left(  t-\delta\right)  ]I_{\left[  \delta,T\right]  }\left(  t\right)
-\alpha(t))Y(t)\}dt]\\
+\mathbb{E[}%
{\textstyle\int_{0}^{T}}
\{(\partial_{u}b)(t)p(t)+(\partial_{u}\sigma)(t)q(t)\}\left(  u\left(
t\right)  -u^{\ast}\left(  t\right)  \right)  dt].
\end{array}
\label{13}%
\end{equation}
Hence, putting%

\begin{equation}%
\begin{array}
[c]{l}%
\zeta(t):=(\partial_{x}b)(t)p(t)+(\partial_{x}\sigma)(t)q(t)\\
+\mathbb{\tilde{E}[}\left(  \partial_{\mu}b\right)  (t-\delta)\tilde{p}\left(
t-\delta\right)  ]I_{\left[  \delta,T\right]  }\left(  t\right)
+\mathbb{\tilde{E}[}\left(  \partial_{\mu}\sigma\right)  (t-\delta)\tilde
{q}\left(  t-\delta\right)  ]I_{\left[  \delta,T\right]  }\left(  t\right)  ,
\end{array}
\label{zetat}%
\end{equation}
and%
\begin{equation}
\zeta:=%
{\textstyle\int_{T-\delta}^{T}}
\mathbb{\tilde{E}[}\left(  \partial_{\mu}b\right)  (t)\tilde{p}\left(
t\right)  +\left(  \partial_{\mu}\sigma\right)  (t)\tilde{q}\left(  t\right)
]dt. \label{zeta}%
\end{equation}
Then, $\left(  \ref{13}\right)  $, takes the form%
\begin{equation}%
\begin{array}
[c]{l}%
\mathbb{E[(}p(T)-\zeta)Y(T)]\\
=\mathbb{E[}\int_{0}^{T}((\partial_{u}b)(t)p(t)+(\partial_{u}\sigma
)(t)q(t)\left(  u\left(  t\right)  -u^{\ast}\left(  t\right)  \right)  dt]\\
+\mathbb{E[}%
{\textstyle\int_{0}^{T}}
(\zeta(t)-\alpha(t))Y(t)dt]\text{.}%
\end{array}
\label{A}%
\end{equation}
We are now able to determine our adjoint process, putting%
\begin{equation}%
\begin{array}
[c]{ll}%
p(T) & :=\zeta+\left(  \partial_{x}g\right)  (X^{\ast}(T),P_{X^{\ast}%
(T)})+\mathbb{\tilde{E}[}\left(  \partial_{\mu}g\right)  (\tilde{X}^{\ast
}(T),P_{X^{\ast}(T)},X^{\ast}(T))]\\
& +%
{\textstyle\int_{T-\delta}^{T}}
\mathbb{\tilde{E}[}\left(  \partial_{\mu}l\right)  (t,\tilde{X}^{\ast
}(t),P_{X^{\ast}(T)},X^{\ast}(T),\tilde{u}^{\ast}(t))]dt,
\end{array}
\label{C}%
\end{equation}
and%
\begin{equation}%
\begin{array}
[c]{ll}%
\alpha(t) & :=\zeta(t)+\left(  \partial_{x}l\right)  (t)+\mathbb{\tilde{E}%
[}\left(  \partial_{\mu}l\right)  (t-\delta)]I_{\left[  \delta,T\right]  }(t),
\end{array}
\label{D}%
\end{equation}
where we denote by%
\[
\left(  \partial_{x}l\right)  (t):=\left(  \partial_{x}l\right)  (t,X^{\ast
}(t),P_{X^{\ast}(t+\delta)},u^{\ast}(t)),
\]%
\[
\mathbb{\tilde{E}[}\left(  \partial_{\mu}l\right)  (t-\delta)]:=\mathbb{\tilde
{E}[}\left(  \partial_{\mu}l\right)  (t-\delta,\tilde{X}(t-\delta),P_{X^{\ast
}(t)},X^{\ast}(t),\tilde{u}^{\ast}(t-\delta))].
\]
Combining $(\ref{zetat})$, $(\ref{zeta})$ with $(\ref{C})$ and $(\ref{D}),$
then $(\ref{bac})$ takes the following form
\begin{equation}%
\begin{array}
[c]{l}%
dp(t)=-\{(\partial_{x}b)(t)p(t)+(\partial_{x}\sigma)(t)q(t)+\left(
\partial_{x}l\right)  (t)\\
+\mathbb{\tilde{E}[}\left(  \partial_{\mu}b\right)  (t-\delta)\tilde{p}\left(
t-\delta\right)  ]I_{\left[  \delta,T\right]  }\left(  t\right)
+\mathbb{\tilde{E}[}\left(  \partial_{\mu}\sigma\right)  (t-\delta)\tilde
{q}\left(  t-\delta\right)  ]I_{\left[  \delta,T\right]  }\left(  t\right) \\
+\mathbb{\tilde{E}[}\left(  \partial_{\mu}l\right)  (t-\delta)]I_{\left[
\delta,T\right]  }\left(  t\right)  \}dt+q(t)dB(t)\text{, }t\in\left[
0,T\right]  \text{,}%
\end{array}
\label{b1}%
\end{equation}
with terminal condition%
\begin{equation}%
\begin{array}
[c]{ll}%
p(T)= & \left(  \partial_{x}g\right)  (X^{\ast}(T),P_{X^{\ast}(T)}%
)+\mathbb{\tilde{E}[}\left(  \partial_{\mu}g\right)  (\tilde{X}^{\ast
}(T),P_{X^{\ast}(T)},X^{\ast}(T))]\\
& +%
{\textstyle\int_{T-\delta}^{T}}
(\mathbb{\tilde{E}[}\left(  \partial_{\mu}b\right)  (t)\tilde{p}\left(
t\right)  ]+\mathbb{\tilde{E}[}\left(  \partial_{\mu}\sigma\right)
(t)\tilde{q}\left(  t\right)  ]+\mathbb{\tilde{E}[}\left(  \partial_{\mu
}l\right)  (t)])dt.
\end{array}
\label{b2}%
\end{equation}
We suppose that the above BSDE (\ref{b1})-(\ref{b2}) has a unique solution
$(p,q)\in S_{\mathbb{F}}^{2}([0,T])\times L_{\mathbb{F}}^{2}([0,T])$. We will
discuss this BSDE in the next section.\newline

\subsection{Stochastic Maximum principle}

We define now the Hamiltonian $H:[0,T]\times\Omega\times%
\mathbb{R}
\times\mathcal{P}_{2}\left(
\mathbb{R}
\right)  \times U\times%
\mathbb{R}
\times%
\mathbb{R}
\rightarrow%
\mathbb{R}
$, as%
\begin{equation}
H(x,\mu,u,p,q)=l(t,x,\mu,u)+b(t,x,\mu,u)p+\sigma(t,x,\mu,u)q. \label{h}%
\end{equation}

\begin{theorem}
[Maximum principle]Let $u^{\ast}(t)$ be an optimal control and $X^{\ast}(t)$
the corresponding trajectory. Then, we have%
\begin{align*}
\partial_{u}H(t,X^{\ast}(t),P_{X^{\ast}(t+\delta)},u^{\ast}%
(t),p(t),q(t))(u(t)-u^{\ast}(t))  &  \geqslant0,\\
dt\text{ }dP\text{-a.e., for all }u  &  \in\mathcal{U},
\end{align*}
where $(p,q)\in S_{\mathbb{F}}^{2}([0,T])\times L_{\mathbb{F}}^{2}([0,T])$ is
the solution of the adjoint equation (\ref{b1})-(\ref{b2}).
\end{theorem}

\noindent{Proof}\quad From (\ref{A}) and (\ref{B}) with the choice (\ref{C})
and (\ref{D}), we get%
\[
0\leq\mathbb{E[}%
{\textstyle\int_{0}^{T}}
\{\left(  \partial_{u}b\right)  (t)p(t)+\left(  \partial_{u}\sigma\right)
(t)q(t)+\left(  \partial_{u}l\right)  (t)\}(u(t)-u^{\ast}(t))dt],
\]
for all $u\in\mathcal{U}$. Assume for some $u\in\mathcal{U},$%
\[%
\begin{array}
[c]{l}%
\Gamma_{u}:=\{\left(  t,\omega\right)  \in\left[  0,T\right]  \times\Omega|\\
\{\left(  \partial_{u}b\right)  (t)p(t)+\left(  \partial_{u}\sigma\right)
(t)q(t)+\left(  \partial_{u}l\right)  (t)(u(t)-u^{\ast}(t))(\omega)<0\}\text{,
}%
\end{array}
\]
is such that%
\[
\mathbb{E[}%
{\textstyle\int_{0}^{T}}
I_{\Gamma_{u}}(t)dt]>0.
\]
Then, for $\tilde{u}(t):=u(t)I_{\Gamma_{u}}(t)+u^{\ast}(t)I_{\Gamma_{u}^{c}%
}(t),t\in\left[  0,T\right]  ,\tilde{u}\in\mathcal{U}$ is such that%
\[%
\begin{array}
[c]{l}%
0\leq\mathbb{E[}%
{\textstyle\int_{0}^{T}}
\{\left(  \partial_{u}b\right)  (t)p(t)+\left(  \partial_{u}\sigma\right)
(t)q(t)+\left(  \partial_{u}l\right)  (t)p(t)\}\times\\
\times(u(t)-u^{\ast}(t))I_{\Gamma_{u}}(t)dt]<0.
\end{array}
\]
But this is a contradiction and proves that%
\[
\{\left(  \partial_{u}b\right)  (t)p(t)+\left(  \partial_{u}\sigma\right)
(t)q(t)+\left(  \partial_{u}l\right)  (t)(u(t)-u^{\ast}(t))\}\geq0,
\]
$dt$ $dP$-a.e, for all $u\in\mathcal{U}.$ By the definition of $H$ in $\left(
\ref{h}\right)  $, the proof is complete.$\square$

\section{Solvability of the delayed McKean-Vlasov BSDE}

We now study the BSDE which is the adjoint equation to the above control
problem. We consider the BSDE%
\[
\left\{
\begin{array}
[c]{ll}%
dp(t) & =-\alpha(t)dt+q(t)dB(t)\text{, }t\in\left[  0,T\right]  \text{,}\\
p(T), &
\end{array}
\right.
\]
which we have seen due to our computations that it has the form%
\begin{equation}%
\begin{array}
[c]{l}%
dp(t)=-\{(\partial_{x}b)(t,X^{\ast}\left(  t\right)  ,P_{X^{^{\ast}}\left(
t+\delta\right)  },u^{\ast}\left(  t\right)  )p(t)\\
+(\partial_{x}\sigma)(t,X^{\ast}\left(  t\right)  ,P_{X^{^{\ast}}\left(
t+\delta\right)  },u^{\ast}\left(  t\right)  )q(t)+\left(  \partial
_{x}l\right)  (X^{\ast}(t),P_{X^{\ast}(t+\delta)},u^{\ast}(t))\\
+\mathbb{\tilde{E}[}\left(  \partial_{\mu}b\right)  \left(  \tilde{X}^{\ast
}\left(  t-\delta\right)  ,P_{X^{\ast}\left(  t\right)  },X^{\ast}\left(
t\right)  ,\tilde{u}^{\ast}\left(  t-\delta\right)  \right)  \tilde{p}\left(
t-\delta\right)  ]I_{\left[  \delta,T\right]  }\left(  t\right) \\
+\mathbb{\tilde{E}[}\left(  \partial_{\mu}\sigma\right)  \left(  \tilde
{X}^{\ast}\left(  t-\delta\right)  ,P_{X^{\ast}\left(  t\right)  },X^{\ast
}\left(  t\right)  ,\tilde{u}^{\ast}\left(  t-\delta\right)  \right)
\tilde{q}\left(  t-\delta\right)  ]I_{\left[  \delta,T\right]  }\left(
t\right) \\
+\mathbb{\tilde{E}[}\left(  \partial_{\mu}l\right)  \left(  \tilde{X}^{\ast
}\left(  t-\delta\right)  ,P_{X^{\ast}\left(  t\right)  },X^{\ast}\left(
t\right)  ,\tilde{u}^{\ast}\left(  t-\delta\right)  \right)  ]I_{\left[
\delta,T\right]  }\left(  t\right)  \}dt+q(t)dB(t)\text{, }t\in\left[
0,T\right]  \text{,}%
\end{array}
\label{3}%
\end{equation}
with%
\begin{equation}%
\begin{array}
[c]{ll}%
p(T)= & (\partial_{x}g)(X^{\ast}(T),P_{X^{\ast}(T)})+\mathbb{\tilde{E}%
[(}\partial_{\mu}g)(\tilde{X}^{\ast}(T),P_{X^{\ast}(T)},X^{\ast}(T))]\\
& +%
{\textstyle\int_{T-\delta}^{T}}
(\mathbb{\tilde{E}[(}\partial_{\mu}b)(\tilde{X}^{\ast}\left(  t\right)
,P_{X^{\ast}\left(  T\right)  },X^{\ast}\left(  T\right)  ,\tilde{u}^{\ast
}\left(  t\right)  )\tilde{p}\left(  t\right)  ]\\
& +\mathbb{\tilde{E}}\left[  (\partial_{\mu}\sigma)\left(  \tilde{X}^{\ast
}\left(  t\right)  ,P_{X^{\ast}\left(  T\right)  },X^{\ast}\left(  T\right)
,\tilde{u}^{\ast}\left(  t\right)  \right)  \tilde{q}\left(  t\right)  \right]
\\
& +\mathbb{\tilde{E}}\left[  (\partial_{\mu}l)(\tilde{X}^{\ast}(t),P_{X^{\ast
}(T)},X^{\ast}(T),\tilde{u}^{\ast}(t))\right]  )dt.
\end{array}
\label{4}%
\end{equation}
Let us better understand the form of this BSDE: for $\left(  t,\omega
,\tilde{\omega}\right)  \in\left[  0,T\right]  \times\Omega\times\tilde
{\Omega},x_{1},x_{2},x_{3},x_{4}\in%
\mathbb{R}
$, putting%
\begin{align*}
\theta_{t}(\omega,\tilde{\omega},x_{1},x_{2},x_{3},x_{4})  &  :=(\partial
_{x}b)(t,X^{\ast}\left(  t,\omega\right)  ,P_{X^{\ast}\left(  t+\delta\right)
},u^{\ast}\left(  t,\omega\right)  )x_{1}\\
&  +(\partial_{x}\sigma)(t,X^{\ast}\left(  t,\omega\right)  ,P_{X^{\ast
}\left(  t+\delta\right)  },u^{\ast}\left(  t,\omega\right)  )x_{2}\\
&  +(\partial_{x}l)(t,X^{\ast}\left(  t,\omega\right)  ,P_{X^{\ast}\left(
t+\delta\right)  },u^{\ast}\left(  t,\omega\right)  )\\
&  +(\partial_{\mu}b)(t-\delta,\tilde{X}^{\ast}\left(  t-\delta,\omega\right)
,P_{X^{\ast}\left(  t\right)  },X^{\ast}\left(  t,\omega\right)  ,\tilde
{u}^{\ast}\left(  t-\delta,\tilde{\omega}\right)  )x_{3}I_{\left[
\delta,T\right]  }\left(  t\right) \\
&  +(\partial_{\mu}\sigma)(t-\delta,\tilde{X}^{\ast}\left(  t-\delta
,\omega\right)  ,P_{X^{\ast}\left(  t\right)  },X^{\ast}\left(  t,\omega
\right)  ,\tilde{u}^{\ast}\left(  t-\delta,\tilde{\omega}\right)
)x_{4}I_{\left[  \delta,T\right]  }\left(  t\right) \\
&  +(\partial_{\mu}l)(t-\delta,\tilde{X}^{\ast}\left(  t-\delta,\omega\right)
,P_{X^{\ast}\left(  t\right)  },X^{\ast}\left(  t,\omega\right)  ,\tilde
{u}^{\ast}\left(  t-\delta,\tilde{\omega}\right)  )I_{\left[  \delta,T\right]
}\left(  t\right)  \text{,}%
\end{align*}
and in order to describe also the terminal condition of our BSDE, we consider
the coefficient%
\begin{align*}
\vartheta_{t}(\omega,\tilde{\omega},x)  &  :=(\partial_{\mu}b)(t,\tilde
{X}^{\ast}\left(  t,\tilde{\omega}\right)  ,P_{X^{\ast}\left(  T\right)
},X^{\ast}\left(  T,\omega\right)  ,\tilde{u}^{\ast}\left(  t,\tilde{\omega
}\right)  )x_{3}\\
&  +(\partial_{\mu}\sigma)(t,\tilde{X}^{\ast}\left(  t,\tilde{\omega}\right)
,P_{X^{\ast}\left(  T\right)  },X^{\ast}\left(  T,\omega\right)  ,\tilde
{u}^{\ast}\left(  t,\tilde{\omega}\right)  )x_{4}\\
&  +(\partial_{\mu}l)(t,\tilde{X}^{\ast}\left(  t,\omega\right)  ,P_{X^{\ast
}\left(  T\right)  },X^{\ast}\left(  T,\omega\right)  ,\tilde{u}^{\ast}\left(
t,\tilde{\omega}\right)  ),
\end{align*}
We know that%
\begin{align}
\varphi(\tilde{P}_{\zeta})(\omega)  &  =\mathbb{\tilde{E}[}\zeta
](\omega)\label{**}\\
&  =%
{\textstyle\int_{\tilde{\Omega}}}
\zeta(\omega,\tilde{\omega})\tilde{P}(d\tilde{\omega}),\omega\in\Omega\text{,
for }\zeta\in L^{2}\left(  \bar{\Omega},\mathcal{\bar{F}},\bar{P}\right)
\text{.}\nonumber
\end{align}
By using $(\ref{**})$ the BSDE $(\ref{3})-(\ref{4})$ takes the form

\begin{definition}
The BSDE $(p,q)\in S_{\mathbb{F}}^{2}([0,T])\times L_{\mathbb{F}}^{2}([0,T])$
is defined by%
\begin{equation}
\left\{
\begin{array}
[c]{ll}%
dp\left(  t\right)  & =-\varphi(\tilde{P}_{\theta_{t}(p\left(  t\right)
,q\left(  t\right)  ,\tilde{p}\left(  t-\delta\right)  ,\tilde{q}\left(
t-\delta\right)  )})dt+q\left(  t\right)  dB\left(  t\right)  \text{, }%
t\in\left[  0,T\right]  \text{,}\\
p\left(  T\right)  & =\zeta+%
{\textstyle\int_{T-\delta}^{T}}
\varphi(\tilde{P}_{\vartheta_{t}(\tilde{p}\left(  t\right)  ,\tilde{q}\left(
t\right)  )})dt,
\end{array}
\right.  \label{dbsde}%
\end{equation}
where%
\[
\zeta:=(\partial_{x}g)(X^{\ast}\left(  T\right)  ,P_{X^{\ast}\left(  T\right)
})+\mathbb{\tilde{E}[(}\partial_{\mu}g)(\tilde{X}^{\ast}\left(  T\right)
,P_{X^{\ast}\left(  T\right)  },X^{\ast}\left(  T\right)  )]\text{,}%
\]
we see that $\zeta\in L^{2}\left(  \Omega,\mathcal{F},P\right)  $.
\end{definition}

\begin{remark}
We call the BSDE $(\ref{dbsde})$, delayed BSDE because the driver at time $t$
depend on both the solution at time $t$ and on its previous value, i.e. the
solution at time $t-\delta$.
\end{remark}

Note that $\theta$ satisfies the following:\newline

\textbf{Assumptions (H.5)}

\begin{enumerate}
\item $\theta:\left[  0,T\right]  \times\Omega\times\tilde{\Omega}\times%
\mathbb{R}
^{4}\rightarrow%
\mathbb{R}
$ is jointly measurable,

\item $\theta_{t}\left(  \cdot,\cdot,x\right)  $ is $\mathcal{F}_{t}%
\otimes\mathcal{\tilde{F}}_{T}$-progressively measurable, for all $x\in%
\mathbb{R}
^{4}$,

\item for all $x,x^{\prime}\in%
\mathbb{R}
^{4}$,%
\[
\left\vert \theta_{t}\left(  \omega,\tilde{\omega},x\right)  -\theta
_{t}\left(  \omega,\tilde{\omega},x^{\prime}\right)  \right\vert \leq
C\left\vert x-x^{\prime}\right\vert \text{, }dtP(d\omega)\tilde{P}%
(d\tilde{\omega})\text{-a.e.}%
\]

\end{enumerate}

Similarly, $\vartheta$ is assumed to satisfy the following:\newline

\textbf{Assumptions (H.6)}

\begin{enumerate}
\item $\vartheta:\left[  T-\delta,T\right]  \times\Omega\times\tilde{\Omega
}\times%
\mathbb{R}
^{2}\rightarrow%
\mathbb{R}
$ is jointly measurable,

\item $\vartheta\left(  \cdot,\cdot,x\right)  $ is $\mathcal{F}_{T}%
\otimes\mathcal{\tilde{F}}_{T}$-measurable, for all $\left(  t,x\right)
\in\left[  T-\delta,T\right]  \times%
\mathbb{R}
^{2}$,

\item $\left\vert \vartheta_{t}\left(  \omega,\tilde{\omega},0\right)
\right\vert \leq C,dtP\left(  d\omega\right)  \tilde{P}\left(  d\tilde{\omega
}\right)  $-a.e, for some constant $C>0$,

\item $|\vartheta_{t}\left(  \omega,\tilde{\omega},x\right)  -\vartheta
_{t}(\omega^{^{\prime}},\tilde{\omega}^{\prime},x^{\prime})|\leq C\left\vert
x-x^{\prime}\right\vert $, for all $x$, $x^{\prime}\in%
\mathbb{R}
^{2}$, $dtP(d\omega)\tilde{P}(d\tilde{\omega})$-a.e.\newline
\end{enumerate}

However, the function $\varphi:\mathcal{P}_{2}(%
\mathbb{R}
)\rightarrow%
\mathbb{R}
$ in a delayed BSDE $(\ref{3})-(\ref{4})$ is Lipschitz continuous.
Consequently, we have the following more general form for our BSDE.\newline

We consider arbitrary $\theta,\vartheta,\varphi,\psi,\zeta$ with $\theta$
satisfying the assumption (H.5), $\vartheta$ satisfying (H.6), $\varphi
,\psi:\mathcal{P}_{2}(%
\mathbb{R}
)\rightarrow%
\mathbb{R}
$ being Lipschitz and $\zeta\in L^{2}\left(  \Omega,\mathcal{F},P\right)  $,
and we study the delayed BSDE,%

\begin{equation}
\left\{
\begin{array}
[c]{ll}%
dp\left(  t\right)  & =-\varphi(\tilde{P}_{\theta_{t}(p\left(  t\right)
,q\left(  t\right)  ,\tilde{p}\left(  t-\delta\right)  ,\tilde{q}\left(
t-\delta\right)  )})dt+q\left(  t\right)  dB\left(  t\right)  \text{, }%
t\in\left[  0,T\right]  \text{,}\\
p\left(  T\right)  & =\zeta+%
{\textstyle\int_{T-\delta}^{T}}
\psi(\tilde{P}_{\vartheta_{t}(\tilde{p}\left(  t\right)  ,\tilde{q}\left(
t\right)  )})dt\text{.}%
\end{array}
\right.  \label{b}%
\end{equation}

\begin{remark}
The adjoint BSDE we describe it above is a special case of $\left(
\ref{b}\right)  $. Indeed, for the adjoint BSDE we have:%
\begin{align*}
\varphi(\tilde{P}_{\vartheta_{t}})(\omega)  &  =\psi(\tilde{P}_{\vartheta_{t}%
})(\omega)=\mathbb{\tilde{E}[}\vartheta_{t}](\omega)\\
&  =%
{\textstyle\int_{\tilde{\Omega}}}
\vartheta_{t}(\omega,\tilde{\omega})\tilde{P}(d\tilde{\omega}),\omega\in
\Omega\text{, for }\vartheta_{t}\in L^{2}\left(  \bar{\Omega},\mathcal{\bar
{F}},\bar{P}\right)  \text{,}%
\end{align*}

\end{remark}

\begin{definition}
We say that $\left(  p,q\right)  \in S_{\mathbb{F}}^{2}\left(  \left[
0,T\right]  \right)  \times L_{\mathbb{F}}^{2}\left(  \left[  0,T\right]
\right)  $ is a solution of $\left(  \ref{b}\right)  $, if%
\[
\left\{
\begin{array}
[c]{lll}%
p\left(  t\right)  & :=p\left(  0\right)  \text{,} & t\leq0\text{,}\\
q\left(  t\right)  & :=0\text{,} & t\leq0\text{,}%
\end{array}
\right.
\]
and if $\left(  \ref{b}\right)  $ is satisfied.
\end{definition}

\begin{theorem}
Under the above assumptions there is some $\delta_{0}>0$ small enough such
that for all $\delta\in\left(  0,\delta_{0}\right]  $, BSDE $\left(
\ref{b}\right)  $ has a unique solution $\left(  p,q\right)  \in
S_{\mathbb{F}}^{2}\left[  0,T\right]  \times L_{\mathbb{F}}^{2}\left(  \left[
0,T\right]  \right)  $.
\end{theorem}

\noindent{Proof}\quad We embed $S_{\mathbb{F}}^{2}\left[  0,T\right]  \subset%
\mathbb{R}
\times L_{\mathbb{F}}^{2}\left(  \left[  0,T\right]  \right)  $: For $U\in
S_{\mathbb{F}}^{2}\left[  0,T\right]  $ we put $U\left(  t\right)  =U\left(
0\right)  $, $t\in\left[  -\delta,0\right]  $, and we observe that%
\[
(U\left(  t\vee0\right)  _{t\in\left[  -\delta,0\right]  }\equiv(U\left(
0\right)  ,\left(  U\left(  t\right)  \right)  _{t\in\left[  0,T\right]
}))\in%
\mathbb{R}
\times L_{\mathbb{F}}^{2}\left(  \left[  0,T\right]  \right)  \text{.}%
\]
For $V\in L_{\mathbb{F}}^{2}\left[  0,T\right]  $ we use the convention that
$V\left(  t\right)  =0$, $t\leq0$. Let $\left(  U,V\right)  \in H=\left(
\mathbb{R}
\times L_{\mathbb{F}}^{2}\left(  \left[  0,T\right]  \right)  \times
L_{\mathbb{F}}^{2}\left(  \left[  0,T\right]  \right)  \right)  $, and
$\left(  p,q\right)  \in S_{\mathbb{F}}^{2}\left(  \left[  0,T\right]
\right)  \times L_{\mathbb{F}}^{2}\left(  \left[  0,T\right]  \right)  \left(
\subset H\right)  $ the unique solution of the equation%
\[
\left\{
\begin{array}
[c]{l}%
dp\left(  t\right)  =-\varphi(\tilde{P}_{\theta_{t}(U\left(  t\right)
,V\left(  t\right)  ,\tilde{U}\left(  t-\delta\right)  ,\tilde{V}\left(
t-\delta\right)  )})dt+q\left(  t\right)  dB\left(  t\right)  ,t\in\left[
0,T\right]  ,\\
p\left(  t\right)  =p(0),q\left(  t\right)  =0,t\leq0,\\
p\left(  T\right)  =\zeta+%
{\textstyle\int_{T-\delta}^{T}}
\varphi(\tilde{P}_{\vartheta_{t}(\tilde{U}\left(  t\right)  ,\tilde{V}\left(
t\right)  )})dt.
\end{array}
\right.
\]
For this observe that the terminal condition is in $L^{2}\left(
\Omega,\mathcal{F},P\right)  $ and the given coefficient of the BSDE is
$\mathbb{F}$-progressively measurable and square integrable. Let us define
\[
\Phi\left(  U,V\right)  :=\left(  p,q\right)  ,\Phi:H\rightarrow H.
\]
For a suitable $\beta>0$ which will be specified later, we define the norm%
\[
\left\Vert \left(  U,V\right)  \right\Vert _{\beta}:=(\mathbb{E[}U^{2}\left(
0\right)  ]+\mathbb{E[}%
{\textstyle\int_{0}^{T}}
e^{\beta t}(\left\vert U\left(  t\right)  \right\vert ^{2}+\left\vert V\left(
t\right)  \right\vert ^{2})dt])^{\frac{1}{2}},U,V\in H,
\]
which is equivalent to the standard norm $\left\Vert \cdot\right\Vert _{0}$
$\left(  \text{for }\beta=0\right)  $ on $H.$ Note that $\left(  H,\left\Vert
\cdot\right\Vert _{0}\right)  $ is a Banach space, and so is $(H,\left\Vert
\cdot\right\Vert _{\beta})$. We show that for some $\delta_{0}>0$, we have for
all $\delta\in\left(  0,\delta_{0}\right]  $ that%
\[
\Phi:(H,\left\Vert \cdot\right\Vert _{\beta})\rightarrow(H,\left\Vert
\cdot\right\Vert _{\beta})
\]
is a contraction, i.e, there is a unique fixed point $\left(  p,q\right)  \in
H,$ such that $\Phi\left(  p,q\right)  =\left(  p,q\right)  .$ Then $\left(
p,q\right)  $ solves BSDE $\left(  \ref{b}\right)  $ and belongs in particular
to $S_{\mathbb{F}}^{2}\left(  \left[  0,T\right]  \right)  .$ Let $\left(
U^{i},V^{i}\right)  \in H,i=1,2,$ and consider $\left(  p^{i},q^{i}\right)
=\Phi\left(  U^{i},V^{i}\right)  $, i.e.,%
\[
\left\{
\begin{array}
[c]{l}%
dp^{i}\left(  t\right)  =-\varphi(\tilde{P}_{\theta_{t}(U^{i}\left(  t\right)
,V^{i}\left(  t\right)  ,\tilde{U}^{i}\left(  t-\delta\right)  ,\tilde{V}%
^{i}\left(  t-\delta\right)  )})dt+q^{i}\left(  t\right)  dB\left(  t\right)
,t\in\left(  0,T\right)  \\
p^{i}\left(  T\right)  =\zeta+%
{\textstyle\int_{T-\delta}^{T}}
\varphi(\tilde{P}_{\vartheta_{t}(\tilde{U}^{i}\left(  t\right)  ,\tilde{V}%
^{i}\left(  t\right)  )})dt\text{, }i=1,2.
\end{array}
\right.
\]
From Itô's formula applied to $e^{\beta t}\left\vert \bar{p}\left(  t\right)
\right\vert ^{2}$, we obtain%
\[%
\begin{array}
[c]{l}%
\mathbb{E[}|\bar{p}\left(  0\right)  |^{2}]+\mathbb{E[}%
{\textstyle\int_{0}^{T}}
e^{\beta t}(\beta\left\vert \bar{p}\left(  t\right)  \right\vert
^{2}+\left\vert \bar{q}\left(  t\right)  \right\vert ^{2})dt]\\
=\mathbb{E[}e^{\beta T}|\bar{p}\left(  T\right)  |]\\
+2\mathbb{E[}%
{\textstyle\int_{0}^{T}}
e^{\beta t}\bar{p}\left(  t\right)  \{\varphi(\tilde{P}_{\theta_{t}%
(U^{1}\left(  t\right)  ,V^{1}\left(  t\right)  ,\tilde{U}^{1}\left(
t-\delta\right)  ,\tilde{V}^{1}\left(  t-\delta\right)  )}-\varphi(\tilde
{P}_{\theta_{t}(U^{2}\left(  t\right)  ,V^{2}\left(  t\right)  ,\tilde{U}%
^{2}\left(  t-\delta\right)  ,\tilde{V}^{2}\left(  t-\delta\right)  )})\}dt].
\end{array}
\]
Observe that, thanks to the assumptions (H.5) and (H.6),%
\[%
\begin{array}
[c]{l}%
|\varphi(\tilde{P}_{\theta_{t}(U^{1}\left(  t\right)  ,V^{1}\left(  t\right)
,\tilde{U}^{1}\left(  t-\delta\right)  ,\tilde{V}^{1}\left(  t-\delta\right)
)})-\varphi(\tilde{P}_{\theta_{t}(U^{2}\left(  t\right)  ,V^{2}\left(
t\right)  ,\tilde{U}^{2}\left(  t-\delta\right)  ,\tilde{V}^{2}\left(
t-\delta\right)  )})|\\
\leq CW_{2}(\tilde{P}_{\theta_{t}(U^{1}\left(  t\right)  ,V^{1}\left(
t\right)  ,\tilde{U}^{1}\left(  t-\delta\right)  ,\tilde{V}^{1}\left(
t-\delta\right)  )},\tilde{P}_{\theta_{t}(U^{2}\left(  t\right)  ,V^{2}\left(
t\right)  ,\tilde{U}^{2}\left(  t-\delta\right)  ,\tilde{V}^{2}\left(
t-\delta\right)  )})\\
\leq C(\mathbb{\tilde{E}[}|\theta_{t}(U^{1}\left(  t\right)  ,V^{1}\left(
t\right)  ,\tilde{U}^{1}\left(  t-\delta\right)  ,\tilde{V}^{1}\left(
t-\delta\right)  )-\theta_{t}(U^{2}\left(  t\right)  ,V^{2}\left(  t\right)
,\tilde{U}^{2}\left(  t-\delta\right)  ,\tilde{V}^{2}\left(  t-\delta\right)
)|^{2}])^{\frac{1}{2}}\\
\leq C(\mathbb{\tilde{E}[}\left\vert \bar{U}\left(  t\right)  \right\vert
^{2}+\left\vert \bar{V}\left(  t\right)  \right\vert ^{2}+\left\vert \bar
{U}\left(  t-\delta\right)  \right\vert ^{2}+\left\vert \bar{V}\left(
t-\delta\right)  \right\vert ^{2}])^{\frac{1}{2}}\\
\leq C(\left\vert \bar{U}\left(  t\right)  \right\vert +\left\vert \bar
{V}\left(  t\right)  \right\vert +(\mathbb{E[}\left\vert \bar{U}\left(
t-\delta\right)  \right\vert ^{2}+\left\vert \bar{V}\left(  t-\delta\right)
\right\vert ^{2}])^{\frac{1}{2}}).
\end{array}
\]
Hence, for some small $\rho>0,$%
\[%
\begin{array}
[c]{l}%
2\mathbb{E[}%
{\textstyle\int_{0}^{T}}
e^{\beta t}\bar{p}\left(  t\right)  \{\varphi(\tilde{P}_{\theta_{t}%
(U^{1}\left(  t\right)  ,V^{1}\left(  t\right)  ,\tilde{U}^{1}\left(
t-\delta\right)  ,\tilde{V}^{1}\left(  t-\delta\right)  )})-\varphi(\tilde
{P}_{\theta_{t}(U^{2}\left(  t\right)  ,V^{2}\left(  t\right)  ,\tilde{U}%
^{2}\left(  t-\delta\right)  ,\tilde{V}^{2}\left(  t-\delta\right)  )})\}dt]\\
\leq C_{\rho}\mathbb{E[}%
{\textstyle\int_{0}^{T}}
e^{\beta t}\left\vert \bar{p}\left(  t\right)  \right\vert ^{2}dt]+\rho
C\mathbb{E[}%
{\textstyle\int_{0}^{T}}
e^{\beta t}(\left\vert \bar{U}\left(  t\right)  \right\vert ^{2}+\left\vert
\bar{V}\left(  t\right)  \right\vert ^{2})dt]\\
+\rho C\mathbb{E[}%
{\textstyle\int_{0}^{T}}
e^{\beta t}(\left\vert \bar{U}\left(  t-\delta\right)  \right\vert
^{2}+\left\vert \bar{V}\left(  t-\delta\right)  \right\vert ^{2})dt].
\end{array}
\]
Note that%
\[%
\begin{array}
[c]{l}%
\rho C\mathbb{E[}%
{\textstyle\int_{0}^{T}}
e^{\beta t}(\left\vert \bar{U}\left(  t-\delta\right)  \right\vert
^{2}+\left\vert \bar{V}\left(  t-\delta\right)  \right\vert ^{2})dt]\\
\leq\rho Ce^{\beta\delta}\mathbb{E[}%
{\textstyle\int_{0}^{T-\delta}}
e^{\beta t}(\left\vert \bar{U}\left(  t\right)  \right\vert ^{2}+\left\vert
\bar{V}\left(  t\right)  \right\vert ^{2})dt]+\rho Ce^{\beta\delta}%
\mathbb{E[}\left\vert \bar{U}\left(  0\right)  \right\vert ^{2}].
\end{array}
\]
Moreover, recall that, $\bar{V}\left(  t\right)  =0$, $t\leq0.$ On the other
hand,%
\[%
\begin{array}
[c]{l}%
\mathbb{E[}e^{\beta T}\left\vert \bar{p}\left(  T\right)  \right\vert ^{2}]\\
=\mathbb{E[}e^{\beta T}(%
{\textstyle\int_{T-\delta}^{T}}
(\varphi(\tilde{P}_{\vartheta_{t}(\tilde{U}^{1}\left(  t\right)  ,\tilde
{V}^{1}\left(  t\right)  )})-\varphi(\tilde{P}_{\vartheta_{t}(\tilde{U}%
^{2}\left(  t\right)  ,\tilde{V}^{2}\left(  t\right)  )}))dt)^{2}]\\
\leq C\delta e^{\beta\delta}\mathbb{E[}%
{\textstyle\int_{T-\delta}^{T}}
e^{\beta T}(\left\vert \bar{U}\left(  t\right)  \right\vert ^{2}+\left\vert
\bar{V}\left(  t\right)  \right\vert ^{2})dt].
\end{array}
\]
Letting $0<\delta\leq\rho$, we obtain%
\[%
\begin{array}
[c]{l}%
\mathbb{E[}\left\vert \bar{p}\left(  0\right)  \right\vert ^{2}]+\mathbb{E[}%
\int_{0}^{T}e^{\beta t}(\beta\left\vert \bar{p}\left(  t\right)  \right\vert
^{2}+\left\vert \bar{q}\left(  t\right)  \right\vert ^{2})dt]\\
\leq C\rho e^{\beta\delta}\mathbb{E[}\left\vert \bar{U}\left(  0\right)
\right\vert ^{2}]+C\rho\left(  1+e^{\beta\delta}\right)  \mathbb{E[}%
{\textstyle\int_{0}^{T}}
(\left\vert \bar{U}\left(  t\right)  \right\vert ^{2}+\left\vert \bar
{V}\left(  t\right)  \right\vert ^{2})dt]\\
+C\rho\mathbb{E[}%
{\textstyle\int_{0}^{T}}
\left\vert \bar{p}\left(  t\right)  \right\vert ^{2}dt].
\end{array}
\]
We choose now $\rho=\tfrac{1}{8C},\beta=C_{\rho}+1$ and $\delta_{0}\in
(0,\frac{1}{8C}),$ such that $\tfrac{1+e^{\beta\delta_{0}}}{8}\leq\tfrac{1}%
{2}$. Then, for all $\delta\in\left(  0,\delta_{0}\right]  ,$%
\[%
\begin{array}
[c]{l}%
\mathbb{E[}\left\vert \bar{p}\left(  0\right)  \right\vert ^{2}]+\mathbb{E[}%
{\textstyle\int_{0}^{T}}
e^{\beta t}(\left\vert \bar{p}\left(  t\right)  \right\vert ^{2}+\left\vert
\bar{q}\left(  t\right)  \right\vert ^{2})dt]\\
\leq\tfrac{1}{2}(\mathbb{E[}\left\vert \bar{U}\left(  0\right)  \right\vert
^{2}]+\mathbb{E[}%
{\textstyle\int_{0}^{T}}
(\left\vert \bar{U}\left(  t\right)  \right\vert ^{2}+\left\vert \bar
{V}\left(  t\right)  \right\vert ^{2})dt]),
\end{array}
\]
i.e.,%
\[
\left\Vert \Phi\left(  U^{1},V^{1}\right)  -\Phi\left(  U^{2},V^{2}\right)
\right\Vert _{\beta}^{2}\leq\tfrac{1}{2}\left\Vert \left(  U^{1},V^{1}\right)
-\left(  U^{2},V^{2}\right)  \right\Vert _{\beta}^{2},
\]
for all $\left(  U^{1},V^{1}\right)  ,\left(  U^{2},V^{2}\right)  \in H.$ This
completes the proof.$\qquad\square$

\section{Appendix}

Let $\delta_{0}^{^{\prime}}\in\left(  0,\delta_{0}\right]  $. For simplicity,
we suppose that $b\equiv0$ and that $\sigma\left(  t,x,\mu,u\right)
=\sigma\left(  t,\mu,u\right)  $, because the case $\sigma\left(
t,x,u\right)  $ is well studied in the literature. Among the vast literature,
we refer, for example, to \cite{Bensoussan1} and \cite{Bensoussan2}. Notice
that, we have
\[
\left\{
\begin{array}
[c]{ll}%
X^{\theta}\left(  t\right)   & =x+%
{\textstyle\int_{0}^{t}}
\sigma(s,P_{X^{\theta}(s+\delta)},u^{\theta}\left(  s\right)  )dB\left(
s\right)  \text{, }t\in\left[  0,T\right]  \text{,}\\
X^{\ast}\left(  t\right)   & =x+%
{\textstyle\int_{0}^{t}}
\sigma\left(  s,P_{X^{\ast}(s+\delta)},u^{\ast}\left(  s\right)  \right)
dB(s)\text{, }t\in\left[  0,T\right]  \text{,}%
\end{array}
\right.
\]
and that equation $\left(  \ref{y}\right)  $ writes as follows%
\[
\left\{
\begin{array}
[c]{l}%
Y(t)=\int_{0}^{t}(\mathbb{\tilde{E}[}\left(  \partial_{\mu}\sigma\right)
(s,P_{X^{\ast}(s+\delta)},\tilde{X}^{\ast}\left(  s+\delta\right)  ,u^{\ast
}\left(  s\right)  )\tilde{Y}(s+\delta)]\\
\text{ \ \ \ \ \ \ \ \ \ }+\left(  \partial_{u}\sigma\right)  (s,P_{X^{\ast
}(s+\delta)},u^{\ast}\left(  s\right)  )(u\left(  s\right)  -u^{\ast}\left(
s\right)  ))dB\left(  s\right)  \text{, }t\in\left[  0,T\right]  \text{,}\\
Y\left(  0\right)  =0\text{, }Y(t)=Y(T)\text{, }t\geq T\text{.}%
\end{array}
\right.
\]
Then,%
\begin{equation}%
\begin{array}
[c]{l}%
X^{\theta}\left(  t\right)  -X^{\ast}\left(  t\right)  \\
=%
{\textstyle\int_{0}^{t}}
{\textstyle\int_{0}^{1}}
\partial_{\lambda}[\sigma(s,P_{X^{\ast}(s+\delta)+\lambda\left(  X^{\theta
}(s+\delta)-X^{\ast}(s+\delta)\right)  },u^{\ast}\left(  s\right)
+\lambda\left(  u^{\theta}\left(  s\right)  -u^{\ast}\left(  s\right)
\right)  )]d\lambda dB(s)\\
=%
{\textstyle\int_{0}^{t}}
{\textstyle\int_{0}^{1}}
\{\mathbb{\tilde{E}[(}\partial_{\mu}\sigma)(s,P_{X^{\ast}(s+\delta
)+\lambda\left(  X^{\theta}(s+\delta)-X^{\ast}(s+\delta)\right)  },\tilde
{X}^{\ast}(s+\delta)\\
+\lambda(\tilde{X}^{\theta}(s+\delta)-\tilde{X}^{\ast}(s+\delta)),u^{\ast
}\left(  s\right)  +\lambda\left(  u^{\theta}\left(  s\right)  -u^{\ast
}\left(  s\right)  \right)  )\cdot\\
\cdot(\tilde{X}^{\theta}(s+\delta)-\tilde{X}^{\ast}(s+\delta)))]+\left(
\partial_{u}\sigma\right)  (s,P_{X^{\ast}(s+\delta)+\lambda\left(  X^{\theta
}(s+\delta)-X^{\ast}(s+\delta)\right)  },u^{\ast}\left(  s\right)  \\
+\lambda(u^{\theta}\left(  s\right)  -u^{\ast}\left(  s\right)  ))(u^{\theta
}\left(  s\right)  -u^{\ast}\left(  s\right)  )\}d\lambda dB\left(  s\right)
\text{.}%
\end{array}
\label{x}%
\end{equation}
Let us define $\zeta^{\theta}\left(  t\right)  :=\tfrac{X^{\theta}\left(
t\right)  -X^{\ast}\left(  t\right)  }{\theta}$, $\theta\in\left(  0,1\right]
$. Then equation $\left(  \ref{x}\right)  $ becomes%
\begin{equation}%
\begin{array}
[c]{l}%
\zeta^{\theta}\left(  t\right)  =%
{\textstyle\int_{0}^{t}}
{\textstyle\int_{0}^{1}}
\{\mathbb{\tilde{E}[(}\partial_{\mu}\sigma)(s,P_{X^{\ast}(s+\delta
)+\lambda\theta\zeta^{\theta}\left(  s+\delta\right)  },\tilde{X}^{\ast
}(s+\delta)+\lambda\theta\tilde{\zeta}^{\theta}\left(  s+\delta\right)
,u^{\ast}\left(  s\right)  \\
+\lambda\theta(u(s)-u^{\ast}\left(  s\right)  )\tilde{\zeta}^{\theta}%
(s+\delta)]\\
+(\partial_{u}\sigma)(s,P_{X^{\ast}(s+\delta)+\lambda\theta\zeta^{\theta
}\left(  s+\delta\right)  },u^{\ast}\left(  s\right)  +\lambda\theta(u\left(
s\right)  -u^{\ast}\left(  s\right)  ))(u\left(  s\right)  -u^{\ast}\left(
s\right)  )\}d\lambda dB\left(  s\right)  .
\end{array}
\label{xi}%
\end{equation}
Recall that $\partial_{\mu}\sigma,\partial_{u}\sigma$ and $U$ are bounded.
Then, with the argument given above, for some $0<\delta_{0}^{^{\prime\prime}%
}\leqslant\delta_{0}^{^{\prime}}$ small enough, for all $\delta\in
\lbrack0,\delta_{0}^{^{\prime\prime}}],$%
\[
\mathbb{E[}\sup_{t\in\left[  0,T\right]  }\left\vert \zeta^{\theta
}(t)\right\vert ^{2}]\leqslant C\text{, }\theta\in\left(  0,1\right]  \text{.}%
\]
Indeed, for suitably defined $(\mathcal{F}_{t}\mathcal{\otimes\tilde{F}}%
_{T}\mathcal{)}$-adapted $(\alpha^{\theta}(s)),(\beta^{\theta}(s))$ processes,
depending on $u^{\theta},u^{\ast}$; but bounded by a bound independent of
$\theta,$ $u^{\theta}$ and $u^{\ast}$, we have%
\[
\left\{
\begin{array}
[c]{ll}%
\zeta^{\theta}(t) & =%
{\textstyle\int_{0}^{t}}
\mathbb{\tilde{E}[}(\alpha^{\theta}(s)\tilde{\zeta}^{\theta}(s+\delta
)+\beta^{\theta}(s))]dB(s),t\in\left[  0,T\right]  ,\\
\zeta^{\theta}(t) & =\zeta^{\theta}(T),t\geqslant T.
\end{array}
\right.
\]
Then, for $\beta>0$, with $\zeta^{\theta}(t)=0$,
\[%
\begin{array}
[c]{l}%
\mathbb{E[}e^{-\beta T}\left\vert \zeta^{\theta}(T)\right\vert ^{2}%
]+\beta\mathbb{E[}%
{\textstyle\int_{0}^{T}}
e^{-\beta s}\left\vert \zeta^{\theta}(s)\right\vert ^{2}ds]\\
=\mathbb{E[}%
{\textstyle\int_{0}^{T}}
e^{-\beta s}|\mathbb{\tilde{E}[}\alpha^{\theta}(s)\tilde{\zeta}^{\theta
}(s+\delta)+\beta^{\theta}(s)]|^{2}ds]\\
\leq C+C\mathbb{E[}%
{\textstyle\int_{0}^{T}}
e^{-\beta s}\left\vert \zeta^{\theta}(s+\delta)\right\vert ^{2}ds]\\
\leq C+C\exp\{\beta\delta\}(\mathbb{E[}%
{\textstyle\int_{0}^{T}}
e^{-\beta s}\left\vert \zeta^{\theta}(s)\right\vert ^{2}ds]+\delta
\mathbb{E[}e^{-\beta T}\left\vert \zeta^{\theta}(T)\right\vert ^{2}])\\
\leq C+3C\delta\mathbb{E[}e^{-\beta T}\left\vert \zeta^{\theta}(T)\right\vert
^{2}]+3C\mathbb{E[}%
{\textstyle\int_{0}^{T}}
e^{-\beta s}\left\vert \zeta^{\theta}(s)\right\vert ^{2}ds],
\end{array}
\]
where we have supposed that $\beta\delta\leq1$. Letting $\beta\delta\leq1$ we
put $\beta=4C$, and we suppose that $\delta\leq\tfrac{1}{6C}$. Then, for some
constant $C^{^{\prime}}$,%
\[%
\begin{array}
[c]{l}%
\mathbb{E[}e^{-\beta T}\left\vert \zeta^{\theta}(T)\right\vert ^{2}%
]+C\mathbb{E[}%
{\textstyle\int_{0}^{T}}
e^{-\beta s}\left\vert \zeta^{\theta}(s)\right\vert ^{2}ds]\leq C^{^{\prime}%
}\text{, }\theta\in\left(  0,1\right]  \text{,}\\
\mathbb{E[}\underset{t\in\left[  0,T\right]  }{\sup}\left\vert \zeta^{\theta
}(t)\right\vert ^{2}]\leq4\mathbb{E[}%
{\textstyle\int_{0}^{T}}
\left\vert \alpha^{\theta}(s)\zeta^{\theta}(s+\delta)+\beta^{\theta
}(s)\right\vert ^{2}ds]\\
\leq C+C\mathbb{E[}%
{\textstyle\int_{0}^{T}}
\left\vert \zeta^{\theta}(s+\delta)\right\vert ^{2}ds]\\
\leq C+C\mathbb{E[}%
{\textstyle\int_{0}^{T}}
\left\vert \zeta^{\theta}(s)\right\vert ^{2}ds]+C\delta\mathbb{E[}\left\vert
\zeta^{\theta}(T)\right\vert ^{2}]\\
\leq C+C\exp\left\{  \beta T\right\}  +C^{^{\prime}}=C^{^{\prime\prime}%
}\text{, }\theta\in\left(  0,1\right]  \text{.}%
\end{array}
\]
We can write now $\left(  \ref{xi}\right)  $ as%
\[%
\begin{array}
[c]{l}%
\zeta^{\theta}(t)=%
{\textstyle\int_{0}^{t}}
(\mathbb{\tilde{E}[}(\partial_{\mu}\sigma)(s,P_{X^{\ast}(s+\delta)},\tilde
{X}^{\ast}(s+\delta),u^{\ast}(s))\tilde{\zeta}^{\theta}(s+\delta)]\\
\left.  +(\partial_{u}\sigma)(s,P_{X^{\ast}(s+\delta)},u^{\ast}%
(s))(u(s)-u^{\ast}(s))\right)  dB(s)+R^{\theta}(t)\text{,}%
\end{array}
\]
with%
\[%
\begin{array}
[c]{l}%
R^{\theta}(t)=%
{\textstyle\int_{0}^{t}}
{\textstyle\int_{0}^{1}}
\{\mathbb{\tilde{E}[}(\partial_{\mu}\sigma)(s,P_{X^{\ast}(s+\delta
)+\lambda\theta\zeta^{\theta}(s+\delta)},\tilde{X}^{\ast}(s+\delta
)+\lambda\theta\tilde{\zeta}^{\theta}(s+\delta),u^{\ast}(s)\\
+\lambda\theta(u(s)-u^{\ast}(s))-(\partial_{\mu}\sigma)(s,P_{X^{\ast}%
(s+\delta)},\tilde{X}^{\ast}(s+\delta),u^{\ast}(s))\tilde{\zeta}^{\theta
}(s+\delta)]\\
+((\partial_{u}\sigma)(s,P_{X^{\ast}(s+\delta)+\lambda\theta\zeta^{\theta
}(s+\delta)},u^{\ast}(s)+\lambda\theta(u(s)-u^{\ast}(s))\\
-(\partial_{u}\sigma)(s,P_{X^{\ast}(s+\delta)},u^{\ast}(s))(u(s)-u^{\ast
}(s))\}d\lambda dB(s),t\in\left[  0,T\right]  .
\end{array}
\]
We have the following estimate%
\[%
\begin{array}
[c]{l}%
\mathbb{E[}\underset{t\in\left[  0,T\right]  }{\sup}\left\vert R^{\theta
}(t)\right\vert ^{2}]\\
\leq2\mathbb{E\tilde{E}[}%
{\textstyle\int_{0}^{T}}
{\textstyle\int_{0}^{1}}
|\partial_{\mu}\sigma(s,P_{X^{\ast}(s+\delta)+\lambda\theta\zeta^{\theta
}(s+\delta)},\tilde{X}^{\ast}(s+\delta)+\lambda\theta\tilde{\zeta}^{\theta
}(s+\delta),u^{\ast}(s)+\lambda\theta(u(s)\\
-u^{\ast}(s))-(\partial_{\mu}\sigma)(s,P_{X^{\ast}(s+\delta)},\tilde{X}^{\ast
}(s+\delta),u^{\ast}(s))|^{2}|\tilde{\zeta}^{\theta}(s+\delta)|^{2}d\lambda
ds]\\
+2\mathbb{E[}%
{\textstyle\int_{0}^{T}}
{\textstyle\int_{0}^{1}}
|(\partial_{u}\sigma)(s,P_{X^{\ast}(s+\delta)+\lambda\theta\zeta^{\theta
}(s+\delta)},u^{\ast}(s)+\lambda\theta(u(s)-u^{\ast}(s)))\\
-(\partial_{u}\sigma)(s,P_{X^{\ast}(s+\delta)},u^{\ast}(s))(u(s)-u^{\ast
}(s))|^{2}d\lambda ds]\\
=I_{1}^{\theta}+I_{2}^{\theta},
\end{array}
\]
where%
\[%
\begin{array}
[c]{ll}%
I_{1}^{\theta} & =2\mathbb{E[\tilde{E}[}%
{\textstyle\int_{0}^{T}}
{\textstyle\int_{0}^{1}}
|(\partial_{\mu}\sigma)(s,P_{X^{\ast}(s+\delta)+\lambda\theta\zeta^{\theta
}(s+\delta)},\tilde{X}^{\ast}(s+\delta)+\lambda\theta\tilde{\zeta}^{\theta
}(s+\delta),\\
& u^{\ast}(s)+\lambda\theta(u(s)-u^{\ast}(s))-(\partial_{\mu}\sigma
)(s,P_{X^{\ast}(s+\delta)},\tilde{X}^{\ast}(s+\delta),u^{\ast}(s))|^{2}%
|\tilde{\zeta}^{\theta}(s+\delta)|^{2}d\lambda ds]],
\end{array}
\]%
\[%
\begin{array}
[c]{ll}%
I_{2}^{\theta} & =2\mathbb{E[}%
{\textstyle\int_{0}^{T}}
{\textstyle\int_{0}^{1}}
|(\partial_{u}\sigma)(s,P_{X^{\ast}(s+\delta)+\lambda\theta\zeta^{\theta
}(s+\delta)},u^{\ast}(s)+\lambda\theta(u(s)-u^{\ast}(s)))\\
& -(\partial_{u}\sigma)(s,P_{X^{\ast}(s+\delta)},u^{\ast}(s))(u(s)-u^{\ast
}(s))|^{2}d\lambda ds].
\end{array}
\]
Then%
\[%
\begin{array}
[c]{ll}%
I_{1}^{\theta} & \leq2\mathbb{E[\tilde{E}[}\underset{t\in\left[  0,T\right]
}{\sup}|\tilde{\zeta}^{\theta}(t)|^{2}%
{\textstyle\int_{0}^{T}}
{\textstyle\int_{0}^{1}}
|(\partial_{\mu}\sigma)(s,P_{X^{\ast}(s+\delta)+\lambda\theta\zeta^{\theta
}(s+\delta)},\tilde{X}^{\ast}(s+\delta)\\
& +\lambda\theta\tilde{\zeta}^{\theta}(s+\delta),u^{\ast}(s)+\lambda
\theta(u(s)-u^{\ast}(s)))-(\partial_{\mu}\sigma)(s,P_{X^{\ast}(s+\delta
)},X^{\ast}(s+\delta),u^{\ast}(s))|^{2}d\lambda ds]]\text{.}%
\end{array}
\]
As $(\partial_{\mu}\sigma)$ is bounded and
\[
\mathbb{\tilde{E}[}\underset{t\in\left[  0,T\right]  }{\sup}|\tilde{\zeta
}^{\theta}(t)|^{2}]\leq C<\infty\text{, }\theta\in(0,1]\text{.}%
\]
We can conclude that $I_{1}^{\theta}\rightarrow0$, as $\theta\rightarrow0$. On
the other hand, as $\partial_{u}\sigma$ is bounded, it follows from the
bounded convergence theorem that%
\[%
\begin{array}
[c]{ll}%
I_{2}^{\theta} & \leq C\mathbb{E[}%
{\textstyle\int_{0}^{T}}
{\textstyle\int_{0}^{1}}
|(\partial_{\mu}\sigma)(s,P_{X^{\ast}(s+\delta)+\lambda\theta\zeta^{\theta
}(s+\delta)},u^{\ast}(s)+\lambda\theta(u(s)-u^{\ast}(s))\\
& -(\partial_{u}\sigma)(s,P_{X^{\ast}(s+\delta)},u^{\ast}(s))|^{2}d\lambda
ds]\rightarrow0\text{, as }\theta\rightarrow0\text{.}%
\end{array}
\]
Consequently,%
\[
\mathbb{E[}\underset{t\in\left[  0,T\right]  }{\sup}\left\vert R^{\theta
}(t)\right\vert ^{2}]\rightarrow0\text{, as }\theta\rightarrow0\text{.}%
\]
Recalling that%
\[%
\begin{array}
[c]{ll}%
\tilde{\zeta}^{\theta}(t) & =%
{\textstyle\int_{0}^{t}}
(\mathbb{\tilde{E}[}(\partial_{\mu}\sigma)(s,P_{X^{\ast}(s+\delta)},\tilde
{X}^{\ast}(s+\delta),u^{\ast}(s))\tilde{\zeta}^{\theta}(s+\delta)]\\
& +(\partial_{u}\sigma)(s,P_{X^{\ast}(s+\delta)},u^{\ast}(s)(u(s)-u^{\ast
}(s))))dB(s)+R^{\theta}(t)\text{,}%
\end{array}
\]
and the equation satisfied by $Y$, we see that%
\[%
\begin{array}
[c]{ll}%
Y(t)-\zeta^{\theta}(t) & =%
{\textstyle\int_{0}^{t}}
\mathbb{\tilde{E}[}(\partial_{\mu}\sigma)(s,P_{X^{\ast}(s+\delta)},\tilde
{X}^{\ast}(s+\delta),u^{\ast}(s))(\tilde{Y}(s+\delta)-\tilde{\zeta}^{\theta
}(s+\delta)]dB(s)-R^{\theta}(t)\text{.}%
\end{array}
\]
Consequently, for $\delta\leq\delta_{0}^{^{\prime}}$,%
\[%
\begin{array}
[c]{lll}%
\mathbb{E[}\underset{t\in\left[  0,T\right]  }{\sup}\left\vert Y(t)-\zeta
^{\theta}(t)\right\vert ^{2}] & \leq C\mathbb{E[}\underset{t\in\left[
0,T\right]  }{\sup}\left\vert R^{\theta}(t)\right\vert ^{2}]\rightarrow
0\text{,} & \text{as }\theta\rightarrow0\text{.}%
\end{array}
\text{ }%
\]
This completes the proof.\newline\newline

\textbf{Acknowledgement}\textit{ I would like to thank Rainer Buckdahn for
helpful discussions during my visit to Brest as well as Bernt Øksendal and
Frank Proske.}

\end{document}